\RequirePackage{xstring}
\RequirePackage{etoolbox}

\newif\ifzh
\IfBeginWith{en}{zh}{\zhtrue}{\zhfalse}

\ifzh
    \IfBeginWith{article}{book}{
    \def\mydocumentclass{ctexbook}
  }{\IfBeginWith{article}{report}{
    \def\mydocumentclass{ctexrep}
  }{
    \def\mydocumentclass{ctexart}
  }}
  \else
    \def\mydocumentclass{article}
  \fi

\ifzh

\else

\fi

\ifzh
  \documentclass[fontset=fandol,]{\mydocumentclass}
\else
  \documentclass[11pt]{\mydocumentclass}
\fi

\usepackage{bookmark}

\usepackage{fontspec} 
\usepackage[Medium, Scale=MatchLowercase]{juliamono} 

\ifzh
    \ctexset{
        section/format              += \zihao{-4}\raggedright,
        subsection/format           += \zihao{5}\raggedright,
        subsubsection/format        += \zihao{5}\raggedright,
    }
\else

    \usepackage{xeCJK}

\fi

\usepackage{calc}
\usepackage[a4paper,top=2.75cm,bottom=2.25cm,left=2.75cm,right=2.75cm,marginparwidth=1.75cm]{geometry}

\usepackage{subcaption}

\providecommand{\tightlist}{%
  \setlength{\itemsep}{0pt}\setlength{\parskip}{0pt}}\usepackage{longtable,booktabs,array}
\usepackage{calc} 
\usepackage{etoolbox}
\makeatletter
\patchcmd\longtable{\par}{\if@noskipsec\mbox{}\fi\par}{}{}
\makeatother
\IfFileExists{footnotehyper.sty}{\usepackage{footnotehyper}}{\usepackage{footnote}}
\makesavenoteenv{longtable}
\usepackage{graphicx}
\makeatletter
\newsavebox\pandoc@box
\newcommand*\pandocbounded[1]{
  \sbox\pandoc@box{#1}%
  \Gscale@div\@tempa{\textheight}{\dimexpr\ht\pandoc@box+\dp\pandoc@box\relax}%
  \Gscale@div\@tempb{\linewidth}{\wd\pandoc@box}%
  \ifdim\@tempb\p@<\@tempa\p@\let\@tempa\@tempb\fi
  \ifdim\@tempa\p@<\p@\scalebox{\@tempa}{\usebox\pandoc@box}%
  \else\usebox{\pandoc@box}%
  \fi%
}
\def\fps@figure{htbp}
\makeatother

\usepackage[backend=biber,style=alphabetic,autocite=plain]{biblatex}
\usepackage{hyperref} 

\usepackage{float} 
\usepackage{orcidlink}

\usepackage{titling}
\ifzh
    \pretitle{\begin{center} \zihao{-2} \bfseries}
    \posttitle{\par\end{center}\vskip 0.5em}
\fi

\makeatletter
\providecommand{\subtitle}[1]{
    \apptocmd{\@title}{\vskip 0.5em \par {\ifzh \zihao{4} \kaishu \else \large \fi #1 \par}}{}{}
}
\makeatother

\ifzh
    \preauthor{\begin{center}\zihao{-4}\kaishu 
        \lineskip 0.5em%
        \begin{tabular}[t]{c}}
\else
    \preauthor{\begin{center}\begin{tabular}[t]{c}}
\fi

\postauthor{
    \end{tabular}\par\end{center}
            \begin{center}\small
                    \textsuperscript{1}~%
            School of Mathematics And Statistics, Beijing Institute of
Technology%
            \par%
                    \textsuperscript{2}~%
            Department of Mathematics, National University of Singapore%
            \par%
                    \end{center}
    }

\ifzh
    \predate{}
    \postdate{}
\else
    \predate{}
    \postdate{}
\fi

\ifzh
    \renewenvironment{abstract}{\zihao{-5} \textbf{摘要} \kaishu}{}
\else
  \renewenvironment{abstract}{\noindent \small \textbf{Abstract.}}{}
\fi

\usepackage{amsmath} 
\usepackage{amssymb} 

\usepackage{tikz}
\usepackage{tikz-cd}
\usepackage{quiver}

\usepackage{amsthm} 

\newcommand{\mynewthm}[2]{\newtheorem{#1}[counter]{#2}}

\ifzh 

\newtheoremstyle{suncnthm}
{}
{}
{}
{}
{}
{\phantom{.}}
{5pt plus 1pt minus 1pt}
{\bfseries \thmname{#1}\thmnumber{ #2}\thmnote{（#3）}}

\theoremstyle{suncnthm}

\mynewthm{theorem}{定理}
\mynewthm{lemma}{引理}
\mynewthm{corollary}{推论}
\mynewthm{proposition}{命题}
\mynewthm{conjecture}{猜想}
\mynewthm{definition}{定义}
\mynewthm{example}{例}
\mynewthm{exercise}{习题}
\newtheorem*{solution}{解}
\newtheorem*{remark}{注记}
\newtheorem{refsolution}{解}
\newtheorem{refremark}{注记}

\makeatletter
\newenvironment{myproof}[1][\myproofname]{\par
  \pushQED{\qed}%
  \normalfont \topsep6\p@\@plus6\p@\relax
  \trivlist
  \item[\hskip\labelsep
    \bfseries 
    #1\phantom{.}
    ]\ignorespaces
}{%
  \popQED\endtrivlist\@endpefalse
}
\providecommand{\myproofname}{证明}

\makeatother
\renewenvironment{proof}{\begin{myproof}}{\end{myproof}}

\else 

\theoremstyle{definition}
\mynewthm{theorem}{Theorem}
\mynewthm{lemma}{Lemma}
\mynewthm{corollary}{Corollary}
\mynewthm{proposition}{Proposition}
\mynewthm{conjecture}{Conjecture}
\mynewthm{definition}{Definition}
\mynewthm{example}{Example}
\mynewthm{exercise}{Exercise}

\theoremstyle{remark}

\newtheorem*{remark}{Remark}

\fi 

\theoremstyle{definition}
\mynewthm{problem}{Problem} 

\usepackage{qrcode} 

\title{Beyond Nine Gates: Classification of Full-Sided Mahjong Waits
with \(3n\) Ranks}

\author{
Shubin Xue%
\textsuperscript{1}%
\,\orcidlink{0009-0006-5125-3202}%
\\ \texttt{shubinx@qq.com}%
 \and%
Xingyu Zhong%
\textsuperscript{2}%
\,\orcidlink{0009-0004-8876-5696}%
\\ \texttt{sun123zxy@qq.com}%
}

\date{}

\begin{document}

\maketitle

\begin{abstract}
    Mahjong is a traditional tile-based game popular in East and
    Southeast Asia. Combinatorial study of Mahjong hands has been a
    topic of interest in recreational mathematics, with the Nine Gates
    hand being a well-known example. This single-suit Mahjong hand
    \(\{1,1,1,2,3,4,5,6,7,8,9,9,9\}\) consists of \(13\) tiles of a
    single suit ranked \(1\) through \(9\), arranged such that drawing
    any tile from that suit completes a winning hand of a pair and \(4\)
    melds, regardless of the tile's rank. Folklore has circulated that
    Nine Gates is the unique single-suit hand with this \(9\)-sided-wait
    property. While prior works have verified this uniqueness by
    exhaustive computer enumeration across different finite
    configurations, we present a general classification of full-sided
    waits for Mahjong variants with \(3n\) ranks and \(n+1\) melds. For
    every \(n\geq2\), we report exactly \(8\) full-sided waits, of which
    Nine Gates is the unique live one that respects the \(4\)-copy
    bound. Notably, our proof requires no exhaustive enumeration or any
    prior knowledge of the classification result. To address the problem
    in its right generality, a comprehensive combinatorial framework is
    developed, allowing the number of ranks, the number of melds, and
    the upper bound on tile multiplicities to vary independently. The
    reliability of our results is partially supported by a manually
    curated ongoing formalization in Lean 4.
\end{abstract}

\renewcommand{\vec}{\boldsymbol}

\section{Introduction}\label{introduction}

\emph{Mahjong} is a traditional tile game popular in China, Japan, and
other parts of East and Southeast Asia. The game is played with a deck
of tiles ranked \(1\) through \(9\) in three suits, with \(4\) copies of
each tile. A winning hand consists of one pair and \(4\) melds, each
meld being either three identical tiles (a triplet) or three tiles of
consecutive ranks (a sequence). Four players take turns drawing and
discarding tiles, with the goal of completing a winning hand.

Combinatorial study of Mahjong hands has been a topic of interest in
recreational mathematics, with the \emph{Nine Gates} hand being a
well-known example. This single-suit Mahjong hand \[
\{1,1,1,2,3,4,5,6,7,8,9,9,9\}
\] consists of \(13\) tiles of a single suit arranged such that drawing
any tile from that suit completes a winning hand, regardless of the
tile's rank. Such a hand is known as a \emph{nine-sided wait}.
Completing Nine Gates earns an award in the highest scoring category
across many Mahjong variants (see e.g. \autocite{ema_riichi_rules}),
earning it fame among Mahjong players.

Despite its celebrity, the uniqueness of Nine Gates, in the sense that
it is the only \(13\)-tile single-suit hand that wins after drawing any
rank while respecting the \(4\)-copy bound, has largely circulated as
folklore. Existing mathematical and computational studies of Mahjong
hands have largely focused either on algorithmic aspects such as
computing the \emph{deficiency number} (the minimum number of tile
replacements needed to complete a winning hand) and optimal discard
policies \autocite{li_lets_2019,wang_new_2022,yan_fast_2024}, or on
computer-assisted enumeration in the conventional \(13\)-tile game. The
preprint \autocite[Proposition 5.1]{cheng_mathematical_2019} verifies
the uniqueness of Nine Gates by programmatically enumerating all
\(93600\) \(13\)-tile single-suit hands containing at most \(4\) copies
of each rank, and \autocite{stephenson_wait_2022} extends this
computational approach to classify the wait patterns into \(828\)
categories. Independently, \autocite{geng_amusement_2010} reports \(5\)
solutions with a weaker formulation: first draw a \(13\)-tile hand from
the \(36\)-tile suit, then ask whether any tile drawn from the remaining
\(23\) tiles completes a winning hand. His classification is incomplete:
a computer search reveals \(9\) solutions, including \(8\) from
Theorem~\ref{thm-main} and one additional hand
\(\{2,3,3,3,3,4,5,6,7,7,7,7,8\}\).

Beyond the classical \(9\)-rank, \(13\)-tile game, three natural
generalizations arise. First, one may generalize the suit to \(3n\)
ranks while scaling the winning hand proportionally to one pair and
\(n+1\) melds. We found no prior academic work on this generalized
setting, though it has been discussed in online forums. A Zhihu question
collects proposed proofs of the uniqueness result for \(3n\) ranks and
\(n+1\) melds, where Goldenrain's answer
\autocite{goldenrain_zhihu_2026}, in particular, presents a promising
argument as handwritten notes, reducing the problem to a finite
enumeration of small cases. Second, there are other Mahjong variants
that alter the number of melds without changing the ranks: most notably,
\emph{Taiwanese Mahjong} is played with \(9\)-rank, \(16\)-tile hands
that require one pair and \(5\) melds (instead of \(4\)) to win. In this
setting, \autocite[Section 6]{cheng_mathematical_2019} reports \(11\)
nine-sided waits. Third, most of the works above require the \(9\)-sided
wait to be \emph{live}, that is, the completed hands must respect the
\(4\)-copy bound. It's natural to ask how the result changes if this
bound on tile multiplicities is removed.

In this paper, we provide a general combinatorial framework that
addresses all three generalizations in a unified manner, and present a
classification of full-sided waits in the traditional \(3n\)-rank,
\((n+1)\)-meld setting. Notably, our proof requires no exhaustive
enumeration or any prior knowledge of the classification result.
Reliability of our results is partially supported by a manually curated
ongoing formalization in Lean 4\footnote{GitHub repository:
  \url{https://github.com/sun123zxy/ninegates-lean4}}. An online
visualization tool is also available\footnote{\url{https://sun123zxy.github.io/ninegates-viz/}}
for the readers' convenience.

To state the result in its right generality, we study a generalized
single-suit game with \(3n\) ranks where a winning hand consists of one
pair and \(m\) melds; \((n,m)=(3,4)\) recovers the standard \(13\)-tile
Mahjong configuration, while \((n,m)=(3,5)\) recovers \(16\)-tile
Taiwanese Mahjong. We call a hand a \emph{full-sided \((n,m)\)-wait} if
drawing any rank completes it to a winning hand, and a full-sided
\((n,m)\)-wait is \emph{live} if all of its completed hands are
\emph{valid}, i.e., respect the \(4\)-copy upper bound. With this
configuration, the natural generalization of Nine Gates, \[
\{1,1,1,2,3,\ldots,3n-1,3n,3n,3n\}
\] is a live full-sided \((n,n+1)\)-wait. For a more compact
presentation, we introduce the multiplicity vector
\((c_1,\ldots,c_{3n})\) of a hand \(\boldsymbol c\), where \(c_i\)
counts the tiles of rank \(i\). The Nine Gates hand is then
\((3,1^{3n-2},3)\), where \(1^{3n-2}\) denotes \(3n-2\) consecutive
entries equal to \(1\). The main result of this paper is the following
classification theorem:

\begin{theorem}[Classification of full-sided
waits]\protect\hypertarget{thm-main}{}\label{thm-main}

Let \(n\) be a positive integer. All full-sided \((n,n+1)\)-waits are
precisely the following hands. For \(n=1\), there are \(6\) hands, all
of which are live: \[
\begin{matrix}
(3,1,3) & (2,3,2) & (3,3,1) & (1,3,3) & (3,2,2) & (2,2,3)
\end{matrix}
\] For \(n\geq2\), there are \(8\) hands
\(\mathbf{A}_n,\mathbf{B}_n,\mathbf{B}'_n,\mathbf{C}_n,\mathbf{C}'_n,\mathbf{D}_n,\mathbf{D}'_n,\mathbf{E}_n\).
Five are given by \[
\begin{aligned}
\mathbf{A}_n&=(3,1^{3n-2},3) \\
\mathbf{B}_n&=(3,1^{3n-5},4,1,1,0) \\
\mathbf{C}_n&=(1,1,4,1^{3n-5},3,0) \\
\mathbf{D}_n&=(0,1,4,2,1^{3n-6},3,0) \\
\mathbf{E}_n&=(0,1,1,4,1^{3n-8},4,1,1,0) & n\geq3 \\
\mathbf{E}_2&=(0,1,4,4,1,0) & n=2
\end{aligned}
\] The hands \(\mathbf{B}'_n,\mathbf{C}'_n,\mathbf{D}'_n\) are obtained
by reversing \(\mathbf{B}_n,\mathbf{C}_n,\mathbf{D}_n\), respectively;
\(\mathbf{A}_n\) and \(\mathbf{E}_n\) are unchanged by reversal. As a
direct consequence, \(\mathbf{A}_n\) is the unique live full-sided
\((n,n+1)\)-wait for \(n\geq2\).

\end{theorem}

Note that in our formulation, a full-sided wait is not required to be a
valid hand. Nevertheless, all full-sided waits listed in
Theorem~\ref{thm-main} are valid.

The paper is structured as follows. In Section~\ref{sec-formulation}, we
formally model the problem and introduce our notation. In
Section~\ref{sec-meldability}, we study meldability, establishing the
language and tools for the subsequent analysis. In
Section~\ref{sec-draw-pair-relation}, we introduce the draw-pair
relation to study the structure of full-sided waits. In
Section~\ref{sec-classification}, we prove a positivity result and
classify all full-sided \((n,m)\)-waits for \(m \leq n+1\), proving
Theorem~\ref{thm-main}. Finally, Section~\ref{sec-amt} discusses minimal
full-sided waits and their classification problem as a direction of
further research.

\section{Formulating the Problem}\label{sec-formulation}

In this section, we describe the combinatorial model of generalized
single-suit Mahjong hands and formulate the problem of classifying
full-sided waits.

Fix an integer \(n\geq1\), called the \emph{order}. We consider a single
suit with ranks \(1,2,\ldots,3n\). It is convenient to allow negative
multiplicities of ranks, thus we define a hand as follows. A
\emph{virtual \(n\)-hand} is an integer vector
\(\boldsymbol c=(c_1,\ldots,c_{3n})\in\mathbb Z^{3n}\), where \(c_i\) is
the \emph{multiplicity} of rank \(i\). An \emph{\(n\)-hand} is a virtual
\(n\)-hand with \(\boldsymbol c\in\mathbb N^{3n}\), where \(\mathbb N\)
is the set of nonnegative integers; it represents a multiset of ranks,
with \(\boldsymbol c\) as its \emph{multiplicity vector}.

Define a partial order on virtual hands by
\(\boldsymbol d\leq\boldsymbol c\) if and only if \(d_i\leq c_i\) for
every rank \(i\). Thus a virtual hand \(\boldsymbol c\) is a hand if and
only if \(\boldsymbol c\geq\boldsymbol 0\). Properties of virtual hands
are automatically inherited by hands. By abuse of language, sometimes we
attribute properties of hands to virtual hands. In such cases, the
corresponding ``being a hand'' condition is implicitly required.

\emph{Reversing} a virtual hand replaces its multiplicity vector by
\((c_{3n},\ldots,c_1)\). Addition and subtraction of virtual hands are
always taken coordinatewise in \(\mathbb Z^{3n}\). A hand is
\emph{valid} if \(c_i\leq4\) for all \(i\). Unless validity is
explicitly required, multiplicities of a hand are unrestricted in
\(\mathbb N\).

For a virtual \(n\)-hand \(\boldsymbol c\in\mathbb Z^{3n}\), define the
\(\mathbb Z\)-linear functionals \[
|\boldsymbol c|:=\sum_{i=1}^{3n}c_i
\qquad
S(\boldsymbol c):=\sum_{i=1}^{3n}i c_i
\qquad
|\boldsymbol c|_r:=\sum_{\substack{1\leq i\leq3n\\i\equiv r\pmod3}}c_i
\] Here \(|\boldsymbol c|\) is the \emph{size} of \(\boldsymbol c\) (its
total number of tiles when \(\boldsymbol c\) is a hand), and
\(S(\boldsymbol c)\) is its \emph{rank sum}. We call
\(|\boldsymbol c|_r\) the \emph{size} of congruence class \(r\) in
\(\boldsymbol c\): for hands, it counts the tiles whose ranks lie in
that class with multiplicity. The following direct fact will be used
repeatedly.

\begin{proposition}[]\protect\hypertarget{prp-class-size-sum-congruence}{}\label{prp-class-size-sum-congruence}

For any virtual \(n\)-hand \(\boldsymbol c\in\mathbb Z^{3n}\), we have
\[
|\boldsymbol c|=|\boldsymbol c|_0+|\boldsymbol c|_1+|\boldsymbol c|_2
\text{ and }
S(\boldsymbol c)\equiv|\boldsymbol c|_1+2|\boldsymbol c|_2\equiv|\boldsymbol c|_1-|\boldsymbol c|_2\pmod3
\] As a consequence,
\(|\boldsymbol c|_0\equiv|\boldsymbol c|_1\equiv|\boldsymbol c|_2\pmod3\)
if and only if \(|\boldsymbol c|\equiv S(\boldsymbol c)\equiv0\pmod3\).

\end{proposition}

Let \(\boldsymbol e_i\) be the standard basis vector representing a
single tile of rank \(i\). A \emph{pair} is a hand \(2\boldsymbol e_i\).
A \emph{meld} is either a \emph{triplet}
\(\boldsymbol t_i:=3\boldsymbol e_i\), with \(1\leq i\leq3n\), or a
\emph{sequence}
\(\boldsymbol s_i:=\boldsymbol e_i+\boldsymbol e_{i+1}+\boldsymbol e_{i+2}\),
with \(1\leq i\leq3n-2\). Each triplet \(\boldsymbol t_i\) has size
\(3\), rank sum \(3i\), and contributes three tiles to a single
congruence class; each sequence \(\boldsymbol s_i\) has size \(3\), rank
sum \(3i+3\), and contributes one tile to each congruence class.

For an integer \(m\geq0\), an \emph{\((n,m)\)-meldable hand} is an
\(n\)-hand that is a sum of \(m\) melds; it has \(3m\) tiles. An
\(n\)-hand is \emph{meldable} if it is \((n,m)\)-meldable for some
\(m\), including the empty hand when \(m=0\). An \emph{\((n,m)\)-winning
hand} is a sum of one pair and an \((n,m)\)-meldable hand, and
consequently has \(3m+2\) tiles. The pair in a chosen winning
decomposition is called its \emph{winning pair}. Neither the winning
pair nor the meld decomposition is required to be unique.

\emph{Drawing} (a tile of) rank \(i\) changes \(\boldsymbol c\) to
\(\boldsymbol c+\boldsymbol e_i\). For an \(n\)-hand \(\boldsymbol c\)
with \(|\boldsymbol c|=3m+1\), a rank \(i\) is an
\emph{\((n,m)\)-waiting tile} of \(\boldsymbol c\) if
\(\boldsymbol c+\boldsymbol e_i\) is an \((n,m)\)-winning hand. It is a
\emph{live \((n,m)\)-waiting tile} if this completed hand is also valid.
We call \(\boldsymbol c\) an \emph{\((n,m)\)-wait} if it has at least
one waiting tile, and a \emph{live \((n,m)\)-wait} if it has at least
one live waiting tile.

\begin{definition}[]\protect\hypertarget{def-full-sided-wait}{}\label{def-full-sided-wait}

An \(n\)-hand \(\boldsymbol c\) of size \(3m+1\) is a \emph{full-sided
\((n,m)\)-wait} if every rank \(1\leq i\leq3n\) is a waiting tile. It is
a \emph{live full-sided \((n,m)\)-wait} if every rank is a live waiting
tile.

\end{definition}

Equivalently, full-sidedness requires that for every rank \(i\) there be
a rank \(j\) such that
\(\boldsymbol c+\boldsymbol e_i-2\boldsymbol e_j\) is
\((n,m)\)-meldable. The rank \(j\) may depend on \(i\) and need not be
unique. Note that a full-sided wait is live if and only if \(c_i\leq3\)
for every \(i\).

One recovers the usual winning condition by taking \(m=n+1\), and the
standard Nine Gates construction of live full-sided \((n,n+1)\)-waits is
\(\mathbf{A}_n\). The statement of Theorem~\ref{thm-main} amounts to
classifying all full-sided \((n,n+1)\)-waits without assuming validity.

\section{Meldability}\label{sec-meldability}

\subsection{Meld Placements}\label{sec-meld-placement}

By definition, meldability is witnessed by a \emph{meld placement},
which records the numbers of sequences and triplets used in a meld
decomposition. We now make this notion precise.

\begin{definition}[Meld
placement]\protect\hypertarget{def-meld-placement}{}\label{def-meld-placement}

Let \(\boldsymbol c\in\mathbb Z^{3n}\) be a virtual \(n\)-hand. For any
vector \(\sigma=(\sigma_1,\ldots,\sigma_{3n-2})\in\mathbb Z^{3n-2}\), we
extend \(\sigma_i\) by zero outside \(1\leq i\leq3n-2\).

\begin{itemize}
\item
  A \emph{(virtual) sequence placement} of \(\boldsymbol c\) is a vector
  \(\sigma=(\sigma_1,\ldots,\sigma_{3n-2})\) in \(\mathbb N^{3n-2}\)
  (resp. \(\mathbb Z^{3n-2}\)) satisfying
  \(\boldsymbol c=\sum_{i=1}^{3n-2}\sigma_i\boldsymbol s_i\), or
  coordinatewise \(c_i=\sigma_i+\sigma_{i-1}+\sigma_{i-2}\) for every
  \(1\leq i\leq3n\). We say that \(\boldsymbol c\) is \emph{(virtually)
  sequentially meldable} if it has a (virtual) sequence placement.
\item
  A \emph{(virtual) triplet placement} of \(\boldsymbol c\) is a vector
  \(\tau=(\tau_1,\ldots,\tau_{3n})\) in \(\mathbb N^{3n}\) (resp.
  \(\mathbb Z^{3n}\)) satisfying
  \(\boldsymbol c=\sum_{i=1}^{3n}\tau_i\boldsymbol t_i\), or
  coordinatewise \(c_i=3\tau_i\).
\item
  A \emph{(virtual) meld placement} of \(\boldsymbol c\) is a pair
  \((\sigma,\tau)\) in \(\mathbb N^{3n-2}\times\mathbb N^{3n}\) (resp.
  \(\mathbb Z^{3n-2}\times\mathbb Z^{3n}\)) satisfying
  \(\boldsymbol c=\sum_{i=1}^{3n-2}\sigma_i\boldsymbol s_i+\sum_{i=1}^{3n}\tau_i\boldsymbol t_i\),
  or coordinatewise \(c_i=\sigma_i+\sigma_{i-1}+\sigma_{i-2}+3\tau_i\)
  for every \(1\leq i\leq3n\). We say that \(\boldsymbol c\) is
  \emph{(virtually) meldable} if it has a (virtual) meld placement.
\end{itemize}

\end{definition}

A meldable hand may have multiple meld placements,
e.g.~\(3\boldsymbol s_i=\boldsymbol t_i+\boldsymbol t_{i+1}+\boldsymbol t_{i+2}\).
By contrast, sequence placements are unique even over \(\mathbb Z\).

\begin{proposition}[Characterization of virtual sequence
placements]\protect\hypertarget{prp-virtual-seq-placement}{}\label{prp-virtual-seq-placement}

For any virtual \(n\)-hand \(\boldsymbol c\in\mathbb Z^{3n}\) (setting
\(c_k=0\) for \(k\leq0\)), the recurrence
\(c_i=\sigma_i+\sigma_{i-1}+\sigma_{i-2}\) for \(1\leq i\leq3n\) with
\(\sigma_{-1}=\sigma_0=0\) uniquely determines
\(\sigma_1,\ldots,\sigma_{3n}\in\mathbb Z\), given in closed form by \[
\sigma_j=(c_j-c_{j-1})+(c_{j-3}-c_{j-4})+\cdots
\] In particular for the boundary coordinates,
\(\sigma_{3n-1}=|\boldsymbol c|_2-|\boldsymbol c|_1\) and
\(\sigma_{3n}=|\boldsymbol c|_0-|\boldsymbol c|_2\). Consequently:

\begin{enumerate}
\def\labelenumi{\arabic{enumi}.}
\tightlist
\item
  \(\boldsymbol c\) is virtually sequentially meldable if and only if
  \(|\boldsymbol c|_0=|\boldsymbol c|_1=|\boldsymbol c|_2\), in which
  case its virtual sequence placement
  \((\sigma_1,\ldots,\sigma_{3n-2})\) is unique.
\item
  \(\boldsymbol c\) is sequentially meldable if and only if
  \(|\boldsymbol c|_0=|\boldsymbol c|_1=|\boldsymbol c|_2\) and
  \(\sigma_j\geq0\) for all \(1\leq j\leq3n-2\).
\end{enumerate}

\end{proposition}

\begin{proof}
In \(\mathbb Z[[x]]\), let \(\Gamma(x)=\sum_{i=1}^{3n}c_i x^i\) and
\(\Sigma(x)=\sum_{j=1}^{3n}\sigma_j x^j\) be the generating functions of
\(\boldsymbol c\) and \(\sigma\). The recurrence is equivalent to
\(\Gamma(x)\equiv(1+x+x^2)\Sigma(x)\pmod{x^{3n+1}}\), or \[
\Sigma(x)\equiv\frac{\Gamma(x)}{1+x+x^2}=\frac{1-x}{1-x^3}\Gamma(x)=(1-x+x^3-x^4+\cdots)\Gamma(x)\pmod{x^{3n+1}}
\] Comparing coefficients gives uniqueness and the closed form for
\(\sigma_j\).
\end{proof}

\subsection{Canonical Meld Placements}\label{sec-triplet-greedy}

It turns out that the non-uniqueness of virtual meld placements of a
virtual hand is solely due to the relation
\(3 \boldsymbol s_i = \boldsymbol t_i + \boldsymbol t_{i+1} + \boldsymbol t_{i+2}\).
To eliminate this flexibility, we introduce a canonical form of meld
placements that forces the number of sequences to be nonnegative and
minimal. The approach was once informally touched upon in the editorial
of an online programming contest \autocite{atcoder_arc160d_2023}.

\begin{definition}[Canonicity]\protect\hypertarget{def-canonical-meld-placement}{}\label{def-canonical-meld-placement}

A (virtual) meld placement \((\sigma,\tau)\) of a (virtual) \(n\)-hand
\(\boldsymbol c\) is \emph{canonical} if \(0\leq\sigma_i<3\) for every
\(1\leq i\leq3n-2\).

\end{definition}

If a meld placement is known a priori, the canonical placement can be
explicitly obtained by the following \emph{canonicalization} procedure.

\begin{proposition}[Canonicalization]\protect\hypertarget{prp-meld-placement-canonicalization}{}\label{prp-meld-placement-canonicalization}

Let \((\sigma,\tau)\) be a (virtual) meld placement of
\(\boldsymbol c\in\mathbb N^{3n}\) (resp. \(\mathbb Z^{3n}\)). Its
\emph{canonicalization}, defined by \[
\widehat\sigma_i=\sigma_i\bmod3
\qquad
\widehat\tau_i=\tau_i+\left\lfloor\frac{\sigma_i}{3}\right\rfloor+\left\lfloor\frac{\sigma_{i-1}}{3}\right\rfloor+\left\lfloor\frac{\sigma_{i-2}}{3}\right\rfloor
\] is a canonical (virtual) meld placement of \(\boldsymbol c\).

\end{proposition}

\begin{proof}
For each \(1\leq i\leq3n-2\), write
\(\sigma_i=\widehat\sigma_i+3\lfloor\sigma_i/3\rfloor\) and substitute
\(\boldsymbol t_i+\boldsymbol t_{i+1}+\boldsymbol t_{i+2}\) for all
\(\lfloor\sigma_i/3\rfloor\) copies of \(3\boldsymbol s_i\). Collecting
the coefficients of each \(\boldsymbol s_i\) and \(\boldsymbol t_i\)
yields the displayed formulas, so \((\widehat\sigma,\widehat\tau)\) is a
virtual meld placement of \(\boldsymbol c\) with
\(0\leq\widehat\sigma_i<3\). When \(\sigma,\tau\geq\boldsymbol 0\), all
terms in \(\widehat\tau_i\) are nonnegative, so
\((\widehat\sigma,\widehat\tau)\) is also a canonical meld placement.
\end{proof}

Even without a placement known in advance, the coordinate equation
\(c_i=\sigma_i+\sigma_{i-1}+\sigma_{i-2}+3\tau_i\) together with
\(0\leq\sigma_i<3\) uniquely determines \(\sigma_i\) and \(\tau_i\) from
left to right for \emph{any} virtual hand
\(\boldsymbol c\in\mathbb Z^{3n}\): at each rank \(i\), after supplying
\(\sigma_{i-1}+\sigma_{i-2}\) tiles to sequences already begun, the
number of new sequences \(\sigma_i\) is forced modulo \(3\) and the
remaining tiles form \(\tau_i\) triplets.

\begin{definition}[Triplet-greedy
scan]\protect\hypertarget{def-meldable-scan}{}\label{def-meldable-scan}

For any virtual \(n\)-hand \(\boldsymbol c\in\mathbb Z^{3n}\), its
\emph{triplet-greedy scan} defines
\(\sigma_1,\ldots,\sigma_{3n}\in\{0,1,2\}\) and
\(\tau_1,\ldots,\tau_{3n}\in\mathbb Z\) recursively by initializing
\(\sigma_{-1}=\sigma_0=0\) and setting, for \(i=1,\ldots,3n\), \[
\sigma_i=(c_i-\sigma_{i-1}-\sigma_{i-2})\bmod3
\qquad
\tau_i=\left\lfloor\frac{c_i-\sigma_{i-1}-\sigma_{i-2}}{3}\right\rfloor=\frac{c_i-\sigma_i-\sigma_{i-1}-\sigma_{i-2}}{3}
\]

\end{definition}

Because the recurrence
\(c_i\equiv\sigma_i+\sigma_{i-1}+\sigma_{i-2}\pmod3\) coincides modulo
\(3\) with that of Proposition~\ref{prp-virtual-seq-placement}, reducing
\((1+x+x^2)^{-1}=(1-x)^{-2}=\sum_{m\geq0}(m+1)x^m\) modulo \(3\) yields
an explicit closed form for \(\sigma_j\) and the following
characterization of virtual meldability and meldability.

\begin{proposition}[Criterion of (virtual)
meldability]\protect\hypertarget{prp-canonical-meld-placement}{}\label{prp-canonical-meld-placement}

Let \(\boldsymbol c\in\mathbb Z^{3n}\) be a virtual \(n\)-hand with
triplet-greedy scan parameters
\(\sigma_1,\ldots,\sigma_{3n}\in\{0,1,2\}\) and
\(\tau_1,\ldots,\tau_{3n}\in\mathbb Z\). Then for every
\(1\leq j\leq3n\), \[
\sigma_j=\left(\sum_{k=1}^j(j-k+1)c_k\right)\bmod3
\] In particular,
\(\sigma_{3n-1}=(|\boldsymbol c|_2-|\boldsymbol c|_1)\bmod3\) and
\(\sigma_{3n}=(|\boldsymbol c|_0-|\boldsymbol c|_2)\bmod3\). Moreover:

\begin{enumerate}
\def\labelenumi{\arabic{enumi}.}
\item
  The following are equivalent:

  \begin{itemize}
  \tightlist
  \item
    \(\boldsymbol c\) is virtually meldable;
  \item
    \(\boldsymbol c\) has a unique canonical virtual meld placement,
    namely \(((\sigma_1,\ldots,\sigma_{3n-2}),\tau)\);
  \item
    \(\sigma_{3n-1}=\sigma_{3n}=0\);
  \item
    \(|\boldsymbol c|_0\equiv|\boldsymbol c|_1\equiv|\boldsymbol c|_2\pmod3\);
  \item
    \(|\boldsymbol c|\equiv0\pmod3\) and
    \(S(\boldsymbol c)\equiv0\pmod3\).
  \end{itemize}
\item
  \(\boldsymbol c\) is meldable if and only if any of the equivalent
  conditions in (1) holds, and additionally, \(\tau_i\geq0\) for every
  \(1\leq i\leq3n\). \(((\sigma_1,\ldots,\sigma_{3n-2}),\tau)\) is the
  unique canonical meld placement of \(\boldsymbol c\).
\end{enumerate}

\end{proposition}

\begin{proof}
Reducing the generating function identity of
Proposition~\ref{prp-virtual-seq-placement} modulo \(3\) gives the
displayed formula for \(\sigma_j\in\{0,1,2\}\) and the boundary
congruences for \(\sigma_{3n-1},\sigma_{3n}\).

For (1), by Definition~\ref{def-meld-placement} and
Definition~\ref{def-meldable-scan}, any canonical virtual meld placement
of \(\boldsymbol c\) must satisfy
\(c_i=\sigma_i+\sigma_{i-1}+\sigma_{i-2}+3\tau_i\) with
\(\sigma_i\in\{0,1,2\}\) and \(\sigma_{3n-1}=\sigma_{3n}=0\). Induction
on \(i\) forces \(\sigma_i\) and \(\tau_i\) to coincide with the
triplet-greedy scan parameters, so
\(((\sigma_1,\ldots,\sigma_{3n-2}),\tau)\) is the unique candidate, and
it is a virtual meld placement if and only if
\(\sigma_{3n-1}=\sigma_{3n}=0\). By
Proposition~\ref{prp-meld-placement-canonicalization}, \(\boldsymbol c\)
is virtually meldable if and only if it admits a canonical virtual meld
placement. Finally, the boundary formulas show that
\(\sigma_{3n-1}=\sigma_{3n}=0\) is equivalent to
\(|\boldsymbol c|_0\equiv|\boldsymbol c|_1\equiv|\boldsymbol c|_2\pmod3\),
which by Proposition~\ref{prp-class-size-sum-congruence} is equivalent
to \(|\boldsymbol c|\equiv S(\boldsymbol c)\equiv0\pmod3\).

For (2), if \(\sigma_{3n-1}=\sigma_{3n}=0\) and \(\tau_i\geq0\) for all
\(i\), then \(((\sigma_1,\ldots,\sigma_{3n-2}),\tau)\) is a meld
placement of \(\boldsymbol c\). Conversely, if \(\boldsymbol c\) is
meldable, canonicalizing any meld placement via
Proposition~\ref{prp-meld-placement-canonicalization} produces a
canonical meld placement, which by uniqueness in (1) must be
\(((\sigma_1,\ldots,\sigma_{3n-2}),\tau)\) with
\(\sigma_{3n-1}=\sigma_{3n}=0\) and \(\tau_i\geq0\) for all \(i\).
\end{proof}

Because the sequence parameters \(\sigma_i\) depend only on the
multiplicities modulo \(3\), we obtain the following lemma for comparing
two virtual hands.

\begin{proposition}[]\protect\hypertarget{prp-canonical-sequence-congruence}{}\label{prp-canonical-sequence-congruence}

Let \(\boldsymbol c,\boldsymbol d\in\mathbb Z^{3n}\) be virtual
\(n\)-hands with triplet-greedy sequence parameters \(\sigma,\sigma'\),
respectively. For any \(1\leq k\leq3n\), we have \(c_i\equiv d_i\pmod3\)
for every \(1\leq i\leq k\) if and only if \(\sigma_i=\sigma'_i\) for
every \(1\leq i\leq k\).

\end{proposition}

\begin{proof}
This follows from the recurrence in Definition~\ref{def-meldable-scan}
(or the closed form in Proposition~\ref{prp-canonical-meld-placement})
together with \(c_i\equiv\sigma_i+\sigma_{i-1}+\sigma_{i-2}\pmod3\).
\end{proof}

\subsection{Tile Moving}\label{sec-tile-moving}

A recurring operation in our analysis is moving a single tile from rank
\(i\) to rank \(j\), which replaces \(\boldsymbol c\) by
\(\boldsymbol c-\boldsymbol e_i+\boldsymbol e_j\). By
Proposition~\ref{prp-canonical-meld-placement}, if \(\boldsymbol c\) and
\(\boldsymbol c-\boldsymbol e_i+\boldsymbol e_j\) are both virtually
meldable, their difference \(-\boldsymbol e_i+\boldsymbol e_j\) is
virtually meldable and has size \(0\), so
\(j-i=S(-\boldsymbol e_i+\boldsymbol e_j)\equiv0\pmod3\). Thus a tile
move can preserve meldability only within the same congruence class
modulo \(3\).

When \(\boldsymbol c\) is sequentially meldable, uniqueness of virtual
sequence placements (Proposition~\ref{prp-virtual-seq-placement})
immediately determines how the sequence decomposition changes.

\begin{proposition}[Sequential tile
moving]\protect\hypertarget{prp-tile-moving}{}\label{prp-tile-moving}

Let \(\boldsymbol c\) be a sequentially meldable \(n\)-hand, and let
\(i,j\) be distinct ranks with \(c_i\geq1\). Then
\(\boldsymbol c-\boldsymbol e_i+\boldsymbol e_j\) is sequentially
meldable if and only if \(i\equiv j\pmod3\) and the sequence placement
of \(\boldsymbol c\) contains at least one copy of each required
sequence below. Replacing one copy of each required sequence by the
corresponding replacement sequence gives the sequence placement of
\(\boldsymbol c-\boldsymbol e_i+\boldsymbol e_j\).

{\def\LTcaptype{none} 
\begin{longtable}[]{@{}
  >{\raggedright\arraybackslash}p{(\linewidth - 4\tabcolsep) * \real{0.3333}}
  >{\raggedright\arraybackslash}p{(\linewidth - 4\tabcolsep) * \real{0.3333}}
  >{\raggedright\arraybackslash}p{(\linewidth - 4\tabcolsep) * \real{0.3333}}@{}}
\toprule\noalign{}
\begin{minipage}[b]{\linewidth}\raggedright
Direction
\end{minipage} & \begin{minipage}[b]{\linewidth}\raggedright
Required sequences
\end{minipage} & \begin{minipage}[b]{\linewidth}\raggedright
Replacement sequences
\end{minipage} \\
\midrule\noalign{}
\endhead
\bottomrule\noalign{}
\endlastfoot
\(i<j\) &
\(\boldsymbol s_i,\boldsymbol s_{i+3},\ldots,\boldsymbol s_{j-3}\) &
\(\boldsymbol s_{i+1},\boldsymbol s_{i+4},\ldots,\boldsymbol s_{j-2}\) \\
\(i>j\) &
\(\boldsymbol s_{j+1},\boldsymbol s_{j+4},\ldots,\boldsymbol s_{i-2}\) &
\(\boldsymbol s_j,\boldsymbol s_{j+3},\ldots,\boldsymbol s_{i-3}\) \\
\end{longtable}
}

\end{proposition}

\begin{proof}
By Proposition~\ref{prp-virtual-seq-placement},
\(\boldsymbol c-\boldsymbol e_i+\boldsymbol e_j\) is virtually
sequentially meldable if and only if
\(-\boldsymbol e_i+\boldsymbol e_j\) has equal class sizes, which holds
if and only if \(i\equiv j\pmod3\). When \(i<j\) and
\(i\equiv j\pmod3\), the telescoping identity
\(-\boldsymbol e_r+\boldsymbol e_{r+3}=\boldsymbol s_{r+1}-\boldsymbol s_r\)
sums over \(r=i,i+3,\ldots,j-3\) to give the unique virtual sequence
placement of \(-\boldsymbol e_i+\boldsymbol e_j\). Adding this to the
sequence placement \(\sigma\) of \(\boldsymbol c\) gives the unique
virtual sequence placement of
\(\boldsymbol c-\boldsymbol e_i+\boldsymbol e_j\), which is nonnegative
if and only if \(\sigma_r\geq1\) for each \(r=i,i+3,\ldots,j-3\). The
case \(i>j\) follows by reversing the ranks.
\end{proof}

When \(\boldsymbol c\) is meldable rather than sequentially meldable,
triplets in the canonical meld placement can also absorb or supply
sequences via
\(\boldsymbol t_r+\boldsymbol t_{r+1}+\boldsymbol t_{r+2}=3\boldsymbol s_r\).
We now characterize when a tile move preserves general meldability.
Remarkably, any long-range move within a congruence class decomposes
into single steps of length \(3\), and each step updates the canonical
meld placement locally.

\begin{proposition}[Tile
moving]\protect\hypertarget{prp-meldable-tile-moving-step}{}\label{prp-meldable-tile-moving-step}

Let \(\boldsymbol c\) be a meldable \(n\)-hand with canonical meld
placement \((\sigma,\tau)\), and let \(1\leq r\leq3n-3\) with
\(c_r\geq1\). Then \(\boldsymbol c-\boldsymbol e_r+\boldsymbol e_{r+3}\)
is meldable if and only if one of the following \(4\) mutually exclusive
conditions holds; in each case, performing the indicated replacement in
\((\sigma,\tau)\) directly yields the canonical meld placement of
\(\boldsymbol c-\boldsymbol e_r+\boldsymbol e_{r+3}\).

{\def\LTcaptype{none} 
\begin{longtable}[]{@{}
  >{\raggedright\arraybackslash}p{(\linewidth - 2\tabcolsep) * \real{0.5000}}
  >{\raggedright\arraybackslash}p{(\linewidth - 2\tabcolsep) * \real{0.5000}}@{}}
\toprule\noalign{}
\begin{minipage}[b]{\linewidth}\raggedright
Condition
\end{minipage} & \begin{minipage}[b]{\linewidth}\raggedright
Replacement
\end{minipage} \\
\midrule\noalign{}
\endhead
\bottomrule\noalign{}
\endlastfoot
\(\sigma_r\geq1,\ \sigma_{r+1}\leq1\) &
\(\boldsymbol s_r\longmapsto\boldsymbol s_{r+1}\) \\
\(\sigma_r\geq1,\ \sigma_{r+1}=2\) &
\(\boldsymbol s_r+2\boldsymbol s_{r+1}\longmapsto\boldsymbol t_{r+1}+\boldsymbol t_{r+2}+\boldsymbol t_{r+3}\) \\
\(\sigma_r=0,\ \sigma_{r+1}\leq1,\ \tau_r,\tau_{r+1},\tau_{r+2}\geq1\) &
\(\boldsymbol t_r+\boldsymbol t_{r+1}+\boldsymbol t_{r+2}\longmapsto2\boldsymbol s_r+\boldsymbol s_{r+1}\) \\
\(\sigma_r=0,\ \sigma_{r+1}=2,\ \tau_r\geq1\) &
\(\boldsymbol t_r+2\boldsymbol s_{r+1}\longmapsto2\boldsymbol s_r+\boldsymbol t_{r+3}\) \\
\end{longtable}
}

More generally, if \(j>r\) with \(j\equiv r\pmod3\) and
\(\boldsymbol c-\boldsymbol e_r+\boldsymbol e_j\) is meldable, then one
of the \(4\) conditions above holds at \(r\). Consequently,
\(\boldsymbol c-\boldsymbol e_r+\boldsymbol e_{r+3}\) is also meldable,
and moving a tile from \(r\) to \(j\) can be carried out step by step
through \(r+3,r+6,\ldots,j\), remaining meldable at every intermediate
rank.

\end{proposition}

\begin{proof}
For sufficiency when \(j=r+3\), in each row the replacement changes the
hand by \(-\boldsymbol e_r+\boldsymbol e_{r+3}\), keeps all counts
nonnegative under the stated condition, and replaces
\((\sigma_r,\sigma_{r+1})\) by \((\sigma_r-1,\sigma_{r+1}+1)\),
\((\sigma_r-1,0)\), \((2,\sigma_{r+1}+1)\), or \((2,0)\), respectively.
All of these lie in \(\{0,1,2\}^2\), so the resulting meld placement is
already canonical.

Conversely, suppose \(j>r\) with \(j\equiv r\pmod3\) and
\(\boldsymbol c':=\boldsymbol c-\boldsymbol e_r+\boldsymbol e_j\) is
meldable, with canonical meld placement \((\sigma',\tau')\). If
\(\sigma_r\geq1\), the first or second row of the table holds
immediately according as \(\sigma_{r+1}\leq1\) or \(\sigma_{r+1}=2\).
Suppose now that \(\sigma_r=0\). Since \(j\geq r+3\), the hands
\(\boldsymbol c\) and \(\boldsymbol c'\) agree at all ranks
\(k\leq r+2\) except \(c'_r-c_r=-1\), so
Proposition~\ref{prp-canonical-sequence-congruence} gives
\(\sigma'_k=\sigma_k\) for \(k<r\), while the closed form of
Proposition~\ref{prp-canonical-meld-placement} gives \[
\sigma'_r=(\sigma_r-1)\bmod3=2
\qquad
\sigma'_{r+1}=(\sigma_{r+1}+1)\bmod3
\qquad
\sigma'_{r+2}=\sigma_{r+2}
\] Subtracting \(c_i=\sigma_i+\sigma_{i-1}+\sigma_{i-2}+3\tau_i\) from
\(c'_i=\sigma'_i+\sigma'_{i-1}+\sigma'_{i-2}+3\tau'_i\) at
\(i=r,r+1,r+2\) (where \(\sigma'_k-\sigma_k=0\) for \(k<r\) and
\(k=r+2\), while \(\sigma'_r-\sigma_r=2\)) then yields \[
\begin{aligned}
-1 = c'_r-c_r &= 2 + 3(\tau'_r-\tau_r) \\
0 = c'_{r+1}-c_{r+1} &= (\sigma'_{r+1}-\sigma_{r+1}) + 2 + 3(\tau'_{r+1}-\tau_{r+1}) \\
0 = c'_{r+2}-c_{r+2} &= (\sigma'_{r+1}-\sigma_{r+1}) + 2 + 3(\tau'_{r+2}-\tau_{r+2})
\end{aligned}
\] The first equation gives \(\tau'_r-\tau_r=-1\), hence
\(\tau_r\geq1\). If \(\sigma_{r+1}=2\), the fourth row holds; if
\(\sigma_{r+1}\leq1\), then \(\sigma'_{r+1}-\sigma_{r+1}=1\), so the
last two equations give
\(\tau'_{r+1}-\tau_{r+1}=\tau'_{r+2}-\tau_{r+2}=-1\), hence
\(\tau_{r+1},\tau_{r+2}\geq1\), and the third row holds. Induction on
\((j-r)/3\) completes the step-by-step assertion.
\end{proof}

Reversing the ranks gives the symmetric \(4\)-case table for moving a
tile three ranks to the left, from \(r+3\) to \(r\), and shows that any
leftward move from \(i\) to \(j<i\) with \(j\equiv i\pmod3\) likewise
factors into steps.

\begin{proposition}[Interval property of tile
moving]\protect\hypertarget{prp-meldable-tile-moving-interval}{}\label{prp-meldable-tile-moving-interval}

Let \(\boldsymbol c\) be a meldable \(n\)-hand and let \(i\) be a rank
with \(c_i\geq1\). Then the set of ranks \[
\{j\mid j\equiv i\pmod3\text{ and }\boldsymbol c-\boldsymbol e_i+\boldsymbol e_j\text{ is meldable}\}
\] forms a consecutive interval within the congruence class of \(i\).

\end{proposition}

\begin{proof}
By Proposition~\ref{prp-meldable-tile-moving-step}, if \(j>i\) belongs
to this set, then every intermediate rank \(i+3,i+6,\ldots,j\) also
belongs to it; by reversal, the same holds for \(j<i\) and the
intermediate ranks \(i-3,i-6,\ldots,j\). Since \(i\) itself belongs to
the set, the set is a consecutive interval in its congruence class.
\end{proof}

\section{The Draw-Pair Relation}\label{sec-draw-pair-relation}

To study full-sidedness, we keep track of which pairs can be used in a
winning decomposition after each draw. This leads to the following
definition.

\begin{definition}[Draw-pair
relation]\protect\hypertarget{def-draw-pair-relation}{}\label{def-draw-pair-relation}

Let \(\boldsymbol c\) be an \(n\)-hand. Its \emph{draw-pair relation}
relates a rank \(i\) to a rank \(j\) precisely when
\(\boldsymbol c+\boldsymbol e_i-2\boldsymbol e_j\) is meldable; we write
\(i\to_{\boldsymbol c}j\) for this relation, or simply \(i\to j\) when
the hand is clear. This relation can be equivalently presented as a
directed graph, called the \emph{draw-pair graph} of \(\boldsymbol c\),
with the ranks as vertices and an edge from \(i\) to \(j\) whenever
\(i\to j\).

\end{definition}

The draw-pair relation is a convenient way to encode the winning
decompositions of a hand \(\boldsymbol c\) after each draw: \(i \to j\)
iff drawing rank \(i\) admits a winning decomposition with pair
\(2\boldsymbol e_j\). A rank is a waiting tile iff it has an outgoing
edge, so \(\boldsymbol c\) is full-sided iff every vertex has an
outgoing edge. As a winning decomposition may not be unique, a rank
\(i\) may admit an outdegree greater than \(1\). For example, the
\((2,2)\)-winning hand \((2,2,2,2,0,0)\) admits \(2 \boldsymbol e_1\)
and \(2 \boldsymbol e_4\) as its winning pairs. Note that loops are
allowed in the draw-pair graph. A direct observation is that a rank
\(i\) can only point to itself if \(c_i \geq 1\), or another rank \(j\)
if \(c_j \geq 2\).

More restrictions arise when one considers the congruence classes of the
ranks.

\begin{proposition}[]\protect\hypertarget{prp-draw-pair-congruence}{}\label{prp-draw-pair-congruence}

If \(i\to j\) in the draw-pair graph of a hand \(\boldsymbol c\), then
\(i+j\equiv-S(\boldsymbol c)\pmod3\).

\end{proposition}

\begin{proof}
If \(i\to j\), then \(\boldsymbol c+\boldsymbol e_i-2\boldsymbol e_j\)
is meldable, hence by Proposition~\ref{prp-canonical-meld-placement} we
have \(S(\boldsymbol c)+i-2j\equiv0\pmod3\).
\end{proof}

Write \(s:=S(\boldsymbol c)\bmod3\), then the congruence class of a draw
determines the class of its winning pair via the map \(r\mapsto-s-r\) on
the three congruence classes. Note that this map fixes \(s\) and
exchanges the other two classes.

\begin{definition}[]\protect\hypertarget{def-fixed-interchanging-classes}{}\label{def-fixed-interchanging-classes}

For an \(n\)-hand \(\boldsymbol c\), the class
\(s:=S(\boldsymbol c)\bmod3\) is called the \emph{fixed class} of
\(\boldsymbol c\). The other two classes, denoted by \(x\) and \(y\) in
either order, are called the \emph{interchanging classes} of
\(\boldsymbol c\).

\end{definition}

By Proposition~\ref{prp-draw-pair-congruence}, edges starting in the
fixed class end in that class, whereas edges starting in either
interchanging class end in the other. The draw-pair graph therefore
consists of a directed graph on the fixed class and a directed bipartite
graph between the interchanging classes. No edges may occur between the
two parts. In particular, self loops can only occur in the fixed class.

We record a congruence relation on the sizes of the classes for future
reference.

\begin{proposition}[]\protect\hypertarget{prp-class-size-congruence}{}\label{prp-class-size-congruence}

Let \(\boldsymbol c\) be an \(n\)-hand with
\(|\boldsymbol c|\equiv1\pmod3\), fixed class \(s\), and interchanging
classes \(x,y\). Then
\(|\boldsymbol c|_x+1\equiv|\boldsymbol c|_y+1\equiv|\boldsymbol c|_s\pmod3\).
In particular, this holds for every \((n,m)\)-wait.

\end{proposition}

\begin{proof}
Choose any rank \(r\) in the fixed class \(s\). Then the virtual hand
\(\boldsymbol c-\boldsymbol e_r\) satisfies
\(|\boldsymbol c-\boldsymbol e_r|=|\boldsymbol c|-1\equiv0\pmod3\) and
\(S(\boldsymbol c-\boldsymbol e_r)=S(\boldsymbol c)-r\equiv0\pmod3\), so
Proposition~\ref{prp-class-size-sum-congruence} gives
\(|\boldsymbol c|_s-1\equiv|\boldsymbol c|_x\equiv|\boldsymbol c|_y\pmod3\).
\end{proof}

\begin{example}[Draw-pair relation of
\(\mathbf A_n\)]\protect\hypertarget{exm-draw-pair-nine-gates}{}\label{exm-draw-pair-nine-gates}

For \(\mathbf A_n=(3,1^{3n-2},3)\), the fixed class is \(s=2\), and the
interchanging classes are \(0\) and \(1\). Its draw-pair relation
consists exactly of \[
\begin{aligned}
i&\to1 &  \text{ for all } & i\equiv0\pmod3 \\
i&\to3n & \text{ for all } & i\equiv1\pmod3 \\
i&\to i &  \text{ for all } & i\equiv2\pmod3
\end{aligned}
\] Indeed,
\(\mathbf A_n-2\boldsymbol e_1-2\boldsymbol e_{3n}=(1,\ldots,1)\) admits
a sequence decomposition
\(\boldsymbol s_1+\boldsymbol s_4+\cdots+\boldsymbol s_{3n-2}\). By
Proposition~\ref{prp-tile-moving}, moving its tile of rank \(1\) to any
rank \(i\equiv1\pmod3\) preserves sequential meldability. Adding the
triplet \(\boldsymbol t_1\) gives a meld decomposition of
\(\mathbf A_n+\boldsymbol e_i-2\boldsymbol e_{3n}\), proving the second
family of edges. Reversing the ranks proves the first. For
\(i\equiv2\pmod3\), the hand \(\mathbf A_n-\boldsymbol e_i\) consists of
the two endpoint triplets and the intervals of single tiles from \(2\)
to \(i-1\) and from \(i+1\) to \(3n-1\). Both intervals have lengths
divisible by \(3\), so they split into sequences, proving \(i\to i\).
Finally, the congruence condition and the fact that only ranks \(1\) and
\(3n\) have multiplicity at least \(2\) exclude all other edges.

\end{example}

\subsection{Incoming Intervals and Two-Pair
Decompositions}\label{incoming-intervals-and-two-pair-decompositions}

We now study the structure of the draw-pair graph to understand the
structure of full-sided waits. Note that by fixing a winning pair
\(2\boldsymbol e_q\) and varying the draw \(p\), the hands
\(\boldsymbol c+\boldsymbol e_p-2\boldsymbol e_q\) and
\(\boldsymbol c+\boldsymbol e_{p'}-2\boldsymbol e_q\) differ by a single
tile move, so that our tile-moving results from
Section~\ref{sec-tile-moving} apply directly. This idea ultimately leads
to a decomposition theorem for full-sided waits. As a first consequence
of Proposition~\ref{prp-meldable-tile-moving-interval}, the draws that
can use \(2\boldsymbol e_q\) form a consecutive interval in their
congruence class.

\begin{proposition}[Incoming interval
property]\protect\hypertarget{prp-incoming-interval}{}\label{prp-incoming-interval}

Let \(\boldsymbol c\) be an \(n\)-hand and let \(q\) be a rank. Then the
set of incoming neighbors \[
I_q:=\{p\mid p\to_{\boldsymbol c}q\}
\] forms a consecutive interval within the congruence class
\((-S(\boldsymbol c)-q)\bmod3\).

\end{proposition}

\begin{proof}
By Proposition~\ref{prp-draw-pair-congruence}, every \(p\in I_q\)
satisfies \(p\equiv-S(\boldsymbol c)-q\pmod3\). If \(|I_q|\leq1\), the
assertion is immediate. Otherwise, choose \(p\in I_q\) with \(p\neq q\),
and set
\(\boldsymbol d:=\boldsymbol c+\boldsymbol e_p-2\boldsymbol e_q\). Then
\(\boldsymbol d\) is meldable and \(d_p=c_p+1\geq1\). For any rank
\(k\equiv p\pmod3\), we have
\(\boldsymbol c+\boldsymbol e_k-2\boldsymbol e_q=\boldsymbol d-\boldsymbol e_p+\boldsymbol e_k\),
so \(k\in I_q\) if and only if
\(\boldsymbol d-\boldsymbol e_p+\boldsymbol e_k\) is meldable. The
conclusion now follows from
Proposition~\ref{prp-meldable-tile-moving-interval}.
\end{proof}

In Example~\ref{exm-draw-pair-nine-gates}, the bidirectional edge
\(1\leftrightarrow3n\) came from decomposing \(\mathbf A_n\) into two
pairs \(2\boldsymbol e_1+2\boldsymbol e_{3n}\) and a meldable hand. The
same phenomenon occurs for all eight families of full-sided
\((n,n+1)\)-waits in Theorem~\ref{thm-main}. This motivates the
following definition.

\begin{definition}[Two-pair
decomposition]\protect\hypertarget{def-two-pair-decomposition}{}\label{def-two-pair-decomposition}

A \emph{two-pair decomposition} of an \(n\)-hand \(\boldsymbol c\) is a
decomposition
\(\boldsymbol c=2\boldsymbol e_p+2\boldsymbol e_q+\boldsymbol d\), where
\(1\leq p,q\leq3n\) are distinct and \(\boldsymbol d\) is meldable.

\end{definition}

The following proposition shows that each two-pair decomposition exactly
corresponds to a bidirectional edge in the draw-pair relation, just as a
decomposition \(\boldsymbol c=\boldsymbol e_q+\boldsymbol d\) with
\(\boldsymbol d\) meldable corresponds to a self-loop \(q\to q\).

\begin{proposition}[Bidirectional edges and two-pair
decompositions]\protect\hypertarget{prp-bidirectional-two-pair}{}\label{prp-bidirectional-two-pair}

Let \(\boldsymbol c\) be an \(n\)-hand and let \(p,q\) be distinct
ranks. Then \(p\to q\) and \(q\to p\) if and only if \(\boldsymbol c\)
admits a two-pair decomposition
\(\boldsymbol c=2\boldsymbol e_p+2\boldsymbol e_q+\boldsymbol d\).

\end{proposition}

\begin{proof}
If \(\boldsymbol c=2\boldsymbol e_p+2\boldsymbol e_q+\boldsymbol d\)
with \(\boldsymbol d\) meldable, then
\(\boldsymbol c+\boldsymbol e_p-2\boldsymbol e_q=\boldsymbol d+\boldsymbol t_p\)
and
\(\boldsymbol c-2\boldsymbol e_p+\boldsymbol e_q=\boldsymbol d+\boldsymbol t_q\)
are meldable, giving both edges.

Conversely, suppose \(p\to q\) and \(q\to p\). The hands
\(\boldsymbol u:=\boldsymbol c+\boldsymbol e_p-2\boldsymbol e_q\) and
\(\boldsymbol v:=\boldsymbol c-2\boldsymbol e_p+\boldsymbol e_q\) are
meldable and congruent coordinatewise modulo \(3\). By
Proposition~\ref{prp-canonical-sequence-congruence}, their canonical
meld placements have the same sequence part, hence a triplet
\(\boldsymbol t_p\) can be removed from the canonical decomposition of
\(\boldsymbol u\) to give a meld decomposition of
\(\boldsymbol d=\boldsymbol u-\boldsymbol t_p=\boldsymbol c-2\boldsymbol e_p-2\boldsymbol e_q\).
This gives the required two-pair decomposition.
\end{proof}

We next show that in any hand, two consecutive edges \(i\to q\) and
\(q\to j\) force a bidirectional edge between \(q\) and some rank
between \(i\) and \(j\). The key observation is that a single-step tile
move from rank \(r\) to \(r+3\) never destroys existing triplets in the
canonical meld placement unless a triplet at \(r\) itself is consumed.

\begin{proposition}[]\protect\hypertarget{prp-moving-preserves-triplets}{}\label{prp-moving-preserves-triplets}

Let \(\boldsymbol h\) be a meldable \(n\)-hand with canonical meld
placement \((\sigma,\tau)\), and let \(1\leq r\leq3n-3\) with
\(h_r\geq1\). If \(\tau_r=0\) and
\(\boldsymbol h-\boldsymbol e_r+\boldsymbol e_{r+3}\) is meldable with
canonical meld placement \((\sigma',\tau')\), then \(\tau'_k\geq\tau_k\)
for every \(1\leq k\leq3n\).

\end{proposition}

\begin{proof}
Since \(\tau_r=0\), only the first two rows of the table in
Proposition~\ref{prp-meldable-tile-moving-step} can occur. Neither
replacement removes any triplet, so every triplet in \((\sigma,\tau)\)
is preserved in \((\sigma',\tau')\).
\end{proof}

\begin{proposition}[]\protect\hypertarget{prp-consecutive-edges}{}\label{prp-consecutive-edges}

Let \(\boldsymbol c\) be an \(n\)-hand. If \(i\to q\) and \(q\to j\),
then \(i\equiv j\pmod3\), and there exists a rank \(r\) satisfying \[
\min\{i,j\}\leq r\leq\max\{i,j\}
\qquad
r\equiv i\equiv j\pmod3
\] such that \(r\to q\) and \(q\to r\). Here \(r=q\) is allowed, in
which case \(q\to q\) is a self-loop.

\end{proposition}

\begin{proof}
By Proposition~\ref{prp-draw-pair-congruence},
\(i+q\equiv-S(\boldsymbol c)\equiv q+j\pmod3\), so \(i\equiv j\pmod3\).
If \(i=j\), taking \(r=i\) gives the claim. By reversing the ranks if
necessary, we may assume \(i<j\).

For each \(k=i,i+3,\ldots,j\), put
\(\boldsymbol h_k:=\boldsymbol c+\boldsymbol e_q+\boldsymbol e_k\).
Since \(i\to q\), the hand
\(\boldsymbol c+\boldsymbol e_i-2\boldsymbol e_q\) is meldable; adding
\(\boldsymbol t_q\) shows that \(\boldsymbol h_i\) is meldable, and by
Proposition~\ref{prp-meld-placement-canonicalization} its canonical meld
placement contains at least one copy of \(\boldsymbol t_q\). Likewise,
\(q\to j\) implies that
\(\boldsymbol h_j=\boldsymbol c+\boldsymbol e_q-2\boldsymbol e_j+\boldsymbol t_j\)
is meldable and its canonical meld placement contains
\(\boldsymbol t_j\).

Since
\(\boldsymbol h_j=\boldsymbol h_i-\boldsymbol e_i+\boldsymbol e_j\),
Proposition~\ref{prp-meldable-tile-moving-step} shows that every
intermediate hand \(\boldsymbol h_k\) is meldable. Let
\(r\in\{i,i+3,\ldots,j\}\) be the smallest rank whose canonical meld
placement of \(\boldsymbol h_r\) contains \(\boldsymbol t_r\); such a
rank exists because \(\boldsymbol h_j\) has this property. For each step
from \(k\) to \(k+3\) prior to \(r\), the canonical meld placement of
\(\boldsymbol h_k\) has no triplet at \(k\), so
Proposition~\ref{prp-moving-preserves-triplets} ensures that all
existing triplets---and in particular \(\boldsymbol t_q\)---are
preserved. At \(k=r\), the canonical meld placement of
\(\boldsymbol h_r\) therefore contains both \(\boldsymbol t_q\) and
\(\boldsymbol t_r\). Removing \(\boldsymbol t_q\) and
\(\boldsymbol t_r\) respectively shows that \[
\boldsymbol h_r-\boldsymbol t_q=\boldsymbol c+\boldsymbol e_r-2\boldsymbol e_q
\qquad
\boldsymbol h_r-\boldsymbol t_r=\boldsymbol c+\boldsymbol e_q-2\boldsymbol e_r
\] are both meldable, so \(r\to q\) and \(q\to r\).
\end{proof}

Applying Proposition~\ref{prp-consecutive-edges} to a full-sided wait
shows that every rank actually used as a winning pair must either have a
self-loop or lie on a bidirectional edge, yielding the following
two-pair decomposition theorem for all full-sided waits.

\begin{proposition}[Two-pair decomposition of full-sided
waits]\protect\hypertarget{prp-two-pair}{}\label{prp-two-pair}

Let \(\boldsymbol c\) be a full-sided \((n,m)\)-wait.

\begin{enumerate}
\def\labelenumi{\arabic{enumi}.}
\tightlist
\item
  Every rank \(q\) with at least one incoming edge
  (\(I_q\neq\emptyset\)) either has a self-loop \(q\to q\) or lies on a
  bidirectional edge \(q\leftrightarrow r\) for some \(r\neq q\).
\item
  Every rank \(q\) in an interchanging class with \(I_q\neq\emptyset\)
  lies on a bidirectional edge \(q\leftrightarrow r\) with \(r\) in the
  other interchanging class, and therefore participates in a two-pair
  decomposition
  \(\boldsymbol c=2\boldsymbol e_q+2\boldsymbol e_r+\boldsymbol d\),
  where \(\boldsymbol d\) is \((n,m-1)\)-meldable. In particular, every
  full-sided \((n,m)\)-wait \(\boldsymbol c\) admits such a two-pair
  decomposition.
\end{enumerate}

\end{proposition}

\begin{proof}
For (1), suppose \(i\to q\). By full-sidedness, \(q\) also has an
outgoing edge \(q\to j\). Applying
Proposition~\ref{prp-consecutive-edges} to \(i\to q\) and \(q\to j\)
gives a rank \(r\) with \(r\to q\) and \(q\to r\), which is a self-loop
if \(r=q\) and a bidirectional edge if \(r\neq q\).

For (2), if \(q\) lies in an interchanging class, then by
Proposition~\ref{prp-draw-pair-congruence} every rank pointing to or
from \(q\) lies in the other interchanging class, so \(r\neq q\) in (1).
Applying Proposition~\ref{prp-bidirectional-two-pair} to
\(q\leftrightarrow r\) gives
\(\boldsymbol c=2\boldsymbol e_q+2\boldsymbol e_r+\boldsymbol d\) with
\(\boldsymbol d\) meldable, and
\(|\boldsymbol d|=|\boldsymbol c|-4=3(m-1)\) shows that
\(\boldsymbol d\) is \((n,m-1)\)-meldable. Finally, any draw \(i\) in
one interchanging class must have an outgoing edge \(i\to q\) to some
rank \(q\) in the other interchanging class, ensuring that such a rank
\(q\) with \(I_q\neq\emptyset\) exists.
\end{proof}

\begin{remark}
Whereas we will see in Section~\ref{sec-classification} that the
two-pair decomposition is unique for \(m=n+1\) when \(n\geq2\), it need
not be unique when \(m>n+1\). For example,
\(\boldsymbol c=(0,2,2,4,2,3)\) is a full-sided \((2,4)\)-wait: draws
\(1,2,3,4,5,6\) admit winning pairs at ranks \(5,4,6,2,4,6\),
respectively. Both ranks \(2\) and \(5\) in the interchanging class
\(2\) have incoming edges, giving the two distinct decompositions \[
\boldsymbol c = 2\boldsymbol e_2+2\boldsymbol e_4+(2\boldsymbol s_3+\boldsymbol t_6) = 2\boldsymbol e_4+2\boldsymbol e_5+(2\boldsymbol s_2+\boldsymbol t_6)
\]

In the fixed class, it is possible for a rank \(q\) to have a nonempty
incoming interval \(I_q\) with \(q\not\in I_q\). For instance, in the
full-sided \((3,7)\)-wait \((2,2,2,3,3,2,4,1,3)\) with fixed class
\(0\), rank \(9\) has \(I_9 = \{3,6\}\).
\end{remark}

\section{\texorpdfstring{Classification of Full-Sided \((n,m)\)-Waits,
\(m \le n+1\)}{Classification of Full-Sided (n,m)-Waits, m \textbackslash le n+1}}\label{sec-classification}

\subsection{Interior Positivity}\label{interior-positivity}

We now establish a key positivity property of full-sided waits: the
multiplicities of all interior ranks are positive. This begins with the
fact that a rank with multiplicity \(0\) separates the melds on its two
sides, which we already used once in
Example~\ref{exm-draw-pair-nine-gates}. We record it here for later use.

\begin{proposition}[Split at a missing
rank]\protect\hypertarget{prp-meldable-split}{}\label{prp-meldable-split}

Let \(\boldsymbol c\) be an \(n\)-hand, \(1 \leq k \leq 3n\), and
suppose \(c_k=0\). Write
\(\boldsymbol c=\boldsymbol c_L+\boldsymbol c_R\), where
\(\boldsymbol c_L=(c_1,\ldots,c_{k-1},0,\ldots,0)\) and
\(\boldsymbol c_R=(0,\ldots,0,c_{k+1},\ldots,c_{3n})\). Then
\(\boldsymbol c\) is meldable if and only if both \(\boldsymbol c_L\)
and \(\boldsymbol c_R\) are meldable.

\end{proposition}

\begin{proposition}[]\protect\hypertarget{prp-waits-across-missing-rank}{}\label{prp-waits-across-missing-rank}

Let \(\boldsymbol c\) be an \((n,m)\)-wait, \(1\leq k\leq3n\), and
suppose \(c_k=0\). Split at \(k\) as in
Proposition~\ref{prp-meldable-split} to obtain
\(\boldsymbol c_L,\boldsymbol c_R\). Suppose there are two waiting tiles
of \(\boldsymbol c\) of rank \(a\), \(b\) with \(a<k<b\). Then
\(|\boldsymbol c_L|\equiv|\boldsymbol c_R|\equiv2\pmod3\), and every
relation \(i\to_{\boldsymbol c}j\) with \(i\neq k\) satisfies:

\begin{itemize}
\tightlist
\item
  If \(i<k\), then \(j>k\), \(\boldsymbol c_L+\boldsymbol e_i\) is
  meldable, and \(i\equiv-S(\boldsymbol c_L)\pmod3\).
\item
  If \(i>k\), then \(j<k\), \(\boldsymbol c_R+\boldsymbol e_i\) is
  meldable, and \(i\equiv-S(\boldsymbol c_R)\pmod3\).
\end{itemize}

That is, the pair \(2 \boldsymbol e_j\) is always taken from the
opposite side of \(k\) from the draw \(i\), and the congruence class of
\(i\) is determined by the hand on its side.

\end{proposition}

\begin{proof}
Choose \(u,v\) with \(a\to u\) and \(b\to v\). Since \(c_k=0\) and
\(a,b\neq k\), we have \(u,v\neq k\). Applying
Proposition~\ref{prp-meldable-split} to
\(\boldsymbol c+\boldsymbol e_a-2\boldsymbol e_u\) gives \[
(|\boldsymbol c_L|,|\boldsymbol c_R|)\equiv
\begin{cases}
(1,0) & u<k \\
(2,2) & u>k
\end{cases}
\pmod3
\] Applying it to \(\boldsymbol c+\boldsymbol e_b-2\boldsymbol e_v\)
instead gives \[
(|\boldsymbol c_L|,|\boldsymbol c_R|)\equiv
\begin{cases}
(2,2) & v<k \\
(0,1) & v>k
\end{cases}
\pmod3
\] The only common possibility is
\((|\boldsymbol c_L|,|\boldsymbol c_R|)\equiv(2,2)\pmod3\).

Now let \(i\to j\) with \(i<k\). Again \(j\neq k\). If \(j<k\), then
Proposition~\ref{prp-meldable-split} makes \(\boldsymbol c_R\) meldable,
contradicting \(|\boldsymbol c_R|\equiv2\pmod3\). Thus \(j>k\), and
Proposition~\ref{prp-meldable-split} again makes
\(\boldsymbol c_L+\boldsymbol e_i\) meldable. By
Proposition~\ref{prp-canonical-meld-placement},
\(S(\boldsymbol c_L)+i\equiv0\pmod3\). The case \(i>k\) is symmetric.
\end{proof}

\begin{proposition}[Interior
positivity]\protect\hypertarget{prp-interior-positive}{}\label{prp-interior-positive}

Every full-sided \((n,m)\)-wait \(\boldsymbol c\) satisfies \(c_i\geq1\)
for \(2\leq i\leq3n-1\).

\end{proposition}

\begin{proof}
Suppose \(c_k=0\) for some \(2\leq k\leq3n-1\). Full-sidedness supplies
waiting tiles on both sides of \(k\), so
Proposition~\ref{prp-waits-across-missing-rank} applies. Every rank on
either side is a waiting tile, but two consecutive ranks cannot belong
to the same congruence class modulo \(3\). Each side therefore contains
at most one rank. Thus \(k-1\leq1\) and \(3n-k\leq1\), forcing \(n=1\)
and \(k=2\).

In this remaining case, \(c_2=0\) while
Proposition~\ref{prp-waits-across-missing-rank} gives
\(c_1\equiv c_3\equiv2\pmod3\), hence \(S(\boldsymbol c)\equiv2\pmod3\).
By Proposition~\ref{prp-draw-pair-congruence}, drawing rank \(2\) must
use a pair at rank \(2\), which is impossible since \(c_2=0\).
\end{proof}

\begin{definition}[]\protect\hypertarget{def-boundary-interior-multiplicities}{}\label{def-boundary-interior-multiplicities}

We call \(c_2,\ldots,c_{3n-1}\) the \emph{interior multiplicities} and
\(c_1,c_{3n}\) the \emph{boundary multiplicities}. We call congruence
classes \(0\) and \(1\) the \emph{boundary classes}, and class \(2\) the
\emph{interior class}. Note that each boundary class contains a boundary
rank and \(n-1\) interior ranks, while the interior class contains \(n\)
interior ranks.

\end{definition}

Proposition~\ref{prp-interior-positive} shows that every full-sided wait
has positive interior multiplicities, a first step toward the
classification of full-sided waits.

\subsection{\texorpdfstring{Non-Existence for
\(m \le n\)}{Non-Existence for m \textbackslash le n}}\label{non-existence-for-m-le-n}

Proposition~\ref{prp-interior-positive} gives a lower bound on each
congruence class size: \[
|\boldsymbol c|_0\geq n-1 \qquad |\boldsymbol c|_1\geq n-1 \qquad |\boldsymbol c|_2\geq n
\] We shall see that strengthening these bounds rules out any full-sided
\((n,m)\)-waits with \(m\leq n\).

\begin{proposition}[]\protect\hypertarget{prp-incoming-edges-both-sides}{}\label{prp-incoming-edges-both-sides}

Let \(\boldsymbol c\) be an \(n\)-hand. If \(a\to j\) and \(b\to j\)
with \(a<j<b\), then \(c_j\geq3\).

\end{proposition}

\begin{proof}
The relation \(a\to j\) gives \(c_j\geq2\). Suppose \(c_j=2\), and put
\(\boldsymbol h=\boldsymbol c-2\boldsymbol e_j\). Then \(h_j=0\), while
\(\boldsymbol h+\boldsymbol e_a\) and \(\boldsymbol h+\boldsymbol e_b\)
are meldable. Splitting
\(\boldsymbol h=\boldsymbol h_L+\boldsymbol h_R\) at \(j\),
Proposition~\ref{prp-meldable-split} shows that both
\(\boldsymbol h_L+\boldsymbol e_a\) and \(\boldsymbol h_L\) are
meldable. This forces
\(|\boldsymbol h_L|+1\equiv|\boldsymbol h_L|\equiv0\pmod3\), a
contradiction.
\end{proof}

\begin{proposition}[Lower bounds on congruence class
sizes]\protect\hypertarget{prp-class-size-lower-bounds}{}\label{prp-class-size-lower-bounds}

Let \(\boldsymbol c\) be a full-sided \((n,m)\)-wait with fixed class
\(s\) and interchanging classes \(x,y\). Then \[
|\boldsymbol c|_s\geq n \qquad |\boldsymbol c|_x\geq n+1 \qquad |\boldsymbol c|_y\geq n+1
\]

\end{proposition}

\begin{proof}
Suppose \(|\boldsymbol c|_s<n\). Then \(s\) is a boundary class, the
boundary multiplicity in class \(s\) is \(0\), and every interior
multiplicity in that class is \(1\). Drawing that boundary rank
completes a winning hand, forcing a pair in class \(s\) by
Proposition~\ref{prp-draw-pair-congruence}. But every multiplicity in
class \(s\) is now \(1\), a contradiction. Hence
\(|\boldsymbol c|_s\geq n\).

By Proposition~\ref{prp-draw-pair-congruence}, every draw in class \(x\)
must use a winning pair from class \(y\), so some rank
\(q\equiv y\pmod3\) satisfies \(c_q\geq2\). Suppose
\(|\boldsymbol c|_y\leq n\). If \(y\) is the interior class, then
\(|\boldsymbol c|_y = n\) and every interior multiplicity in class \(y\)
is \(1\), contradicting \(c_q\geq2\). Hence \(y\) is a boundary class.
To ensure the existence of \(c_q\geq2\) against
Proposition~\ref{prp-interior-positive} and \(|\boldsymbol c|_y\leq n\),
the boundary multiplicity in class \(y\) must be \(0\), and \(q\) is the
unique rank of multiplicity \(2\) in its class, with all other interior
multiplicities being \(1\). Consequently, as \(\boldsymbol c\) is a
full-sided wait, \(i\to q\) for every \(i\equiv x\pmod3\).

Since \(q\) is an interior rank of boundary class \(y\), we have
\(3\leq q\leq3n-2\). Choose \(a\in\{q-2,q-1\}\) and \(b\in\{q+1,q+2\}\)
with \(a\equiv b\equiv x\pmod3\). Then \(a\to q\) and \(b\to q\) with
\(a<q<b\), so Proposition~\ref{prp-incoming-edges-both-sides} gives
\(c_q\geq3\), contradicting \(c_q=2\). Thus
\(|\boldsymbol c|_y\geq n+1\). \(|\boldsymbol c|_x\geq n+1\) follows by
symmetry.
\end{proof}

\begin{remark}
Each bound is sharp for every \(n\). The fixed-class bound is attained
by \(\mathbf A_n\), where the fixed class is \(2\) and \[
|\mathbf A_n|_0=|\mathbf A_n|_1=n+2 \qquad |\mathbf A_n|_2=n
\] For \(n=1\), the hand \((2,3,2)\) attains both interchanging-class
bounds. For \(n\geq2\), the hand \(\boldsymbol c=(0,3,3,3,1^{3n-5},3)\)
is a full-sided \((n,n+2)\)-wait with fixed class \(2\) and
\(|\boldsymbol c|_1=n+1\). Indeed, its draw-pair relation contains
\(i\to4\) for \(i\equiv0\pmod3\), \(i\to3n\) for \(i\equiv1\pmod3\), and
\(i\to i\) for \(i\equiv2\pmod3\), as can be checked by direct meld
decompositions.
\end{remark}

\begin{proposition}[]\protect\hypertarget{prp-no-m-le-n-full-wait}{}\label{prp-no-m-le-n-full-wait}

For every \(n\geq1\) and \(0\leq m\leq n\), no full-sided \((n,m)\)-wait
exists.

\end{proposition}

\begin{proof}
If \(\boldsymbol c\) is a full-sided \((n,m)\)-wait, then
Proposition~\ref{prp-class-size-lower-bounds} gives \[
3m+1=|\boldsymbol c|=|\boldsymbol c|_s+|\boldsymbol c|_x+|\boldsymbol c|_y\geq n+(n+1)+(n+1)=3n+2
\] Hence \(m>n\).
\end{proof}

\subsection{\texorpdfstring{The Case
\(m = n+1\)}{The Case m = n+1}}\label{the-case-m-n1}

We now specialize to the classical case \(m=n+1\), where a full-sided
\((n,n+1)\)-wait has size \(3n+4\). The bounds of
Proposition~\ref{prp-class-size-lower-bounds} leave little freedom in
the three class sizes. We first determine these sizes and the possible
winning pairs for \(n\geq2\), then use Proposition~\ref{prp-tile-moving}
to recover the eight families. The case \(n=1\) will be treated
separately in the proof of Theorem~\ref{thm-main}.

\begin{proposition}[Sizes of congruence
classes]\protect\hypertarget{prp-class-sizes-n-plus-one}{}\label{prp-class-sizes-n-plus-one}

Let \(n\geq2\), and let \(\boldsymbol c\) be a full-sided
\((n,n+1)\)-wait with fixed class \(s\) and interchanging classes
\(x,y\). Then \(|\boldsymbol c|_s=n\) and
\(|\boldsymbol c|_x=|\boldsymbol c|_y=n+2\).

\end{proposition}

\begin{proof}
By Proposition~\ref{prp-class-size-congruence},
\(|\boldsymbol c|_x+1\equiv|\boldsymbol c|_y+1\equiv|\boldsymbol c|_s\pmod3\).
Together with Proposition~\ref{prp-class-size-lower-bounds} and
\(|\boldsymbol c|=3n+4\), this leaves precisely \[
(|\boldsymbol c|_s,|\boldsymbol c|_x,|\boldsymbol c|_y)=(n,n+2,n+2)
\quad\text{or}\quad(n+2,n+1,n+1)
\] Suppose the second possibility holds. Let \(i\equiv x\pmod3\) be any
rank and \(j\) any rank with \(i\to j\), so \(j\equiv y\pmod3\). In
\(\boldsymbol c+\boldsymbol e_i-2\boldsymbol e_j\), class \(y\) has size
\(n-1\), while classes \(x\) and \(s\) each have size \(n+2\). Every
meld decomposition therefore has at most \(n-1\) sequences and contains
at least one triplet in each of classes \(x\) and \(s\). Thus, in terms
of the original multiplicities, either rank \(i\) has multiplicity at
least \(2\) or another rank in class \(x\) has multiplicity at least
\(3\); some rank in the fixed class must also have multiplicity at least
\(3\). Exchanging \(x\) and \(y\) gives the corresponding requirements
for every draw in class \(y\).

If \(s\neq2\), the interior class is an interchanging class of size
\(n+1\). By Proposition~\ref{prp-interior-positive}, it has one
multiplicity \(2\) and all others \(1\). Since \(n\geq2\), choose a rank
\(i\) of multiplicity \(1\) in this class, then \(c_i<2\) and every
other multiplicity in the class is less than \(3\), contradicting the
preceding requirement if we draw \(i\).

If \(s=2\), the fixed class is the interior class and has size \(n+2\).
We have shown that some rank \(r\) in this class has \(c_r\geq3\). Since
all its \(n\) multiplicities are positive, necessarily \(c_r=3\) and
every other multiplicity in this class is \(1\). As every draw in the
interchanging classes requires a triplet in the fixed class (which we
have shown above), that triplet must be \(\boldsymbol t_r\). The hand
\(\boldsymbol h=\boldsymbol c-\boldsymbol t_r\) therefore still waits on
every rank in the interchanging classes, but \(h_r=0\). The ranks
\(r-1\) and \(r+1\) are waiting tiles on opposite sides of \(r\).
Moreover, since \(3n\geq6\), at least one of \(r-2\) and \(r+2\) is a
valid rank. Together with its adjacent rank \(r-1\) or \(r+1\), it gives
two waiting tiles in different congruence classes on the same side of
\(r\), contradicting Proposition~\ref{prp-waits-across-missing-rank}.
\end{proof}

\begin{remark}
The assumption \(n\geq2\) is necessary: the full-sided \((1,2)\)-wait
\((2,3,2)\) in Theorem~\ref{thm-main} has
\((|\boldsymbol c|_s,|\boldsymbol c|_x,|\boldsymbol c|_y)=(3,2,2)\).
\end{remark}

\begin{proposition}[Two
pillars]\protect\hypertarget{prp-unique-pair-ranks}{}\label{prp-unique-pair-ranks}

With \(\boldsymbol c, x,y,s\) and the hypotheses of
Proposition~\ref{prp-class-sizes-n-plus-one}:

\begin{enumerate}
\def\labelenumi{\arabic{enumi}.}
\tightlist
\item
  There are two unique ranks \(p\equiv x \pmod 3\) and
  \(q\equiv y\pmod3\) with \(c_p,c_q\geq2\). In fact, \(c_p,c_q\geq3\).
\item
  In the fixed class \(s\), there is at most one rank of multiplicity
  \(2\), and all other multiplicities are less than \(2\).
\item
  For every \(i\equiv s\pmod3\), each winning decomposition of
  \(\boldsymbol c+\boldsymbol e_i\) consists of a pair in class \(s\),
  the two triplets \(\boldsymbol t_p,\boldsymbol t_q\), and \(n-1\)
  sequences.
\end{enumerate}

\end{proposition}

\begin{proof}
For (2), it follows directly from \(|\boldsymbol c|_s=n\) and
Proposition~\ref{prp-interior-positive}.

For (1) and (3), consider a draw in the fixed class. After removing its
winning pair (which also lies in this class), the class has size
\(n-1\), so there are at most \(n-1\) sequences. Each interchanging
class has size \(n+2\) and must therefore contribute at least a triplet.
As both the draw and the winning pair are taken from the fixed class,
some rank in each interchanging class has multiplicity at least \(3\) in
the original hand. This proves the existence of \(p\) and \(q\) such
that \(c_p,c_q\geq3\). By \(|\boldsymbol c|_x=|\boldsymbol c|_y=n+2\)
and Proposition~\ref{prp-interior-positive}, no other rank in either
interchanging class may have multiplicity greater than \(2\), and
\(c_p, c_q \leq 4\). To summarize, we have shown that \(p\) and \(q\)
are unique in their respective interchanging classes with multiplicity
at least \(3\), partially proving (1). Also, \(\boldsymbol t_p\) and
\(\boldsymbol t_q\) must appear exactly once in every winning
decomposition of a fixed-class draw, squeezing the remaining sequences
to obtain the maximum number \(n-1\). This proves (3).

For the strengthened uniqueness assertion \(c_p,c_q\geq2\) in (1),
remove \(\boldsymbol t_p\) and \(\boldsymbol t_q\) from
\(\boldsymbol c\) to obtain
\(\boldsymbol h=\boldsymbol c-\boldsymbol t_p-\boldsymbol t_q\). Every
rank in class \(s\) remains a waiting tile of \(\boldsymbol h\) with its
winning pair in class \(s\). Suppose class \(x\) contains another rank
of multiplicity at least \(2\). The bound \(|\boldsymbol c|_x=n+2\) and
Proposition~\ref{prp-interior-positive} force \(x\) to be a boundary
class, \(p\) to be interior, and \(c_p=3\): a triplet and another pair
require at least \(n+3\) tiles in the interior class or when the triplet
is at the boundary, and raising \(c_p\) above \(3\) also requires at
least \(n+3\) tiles. Thus \(h_p=0\). Since \(3\leq p\leq3n-2\), there
are ranks \(a<p<b\) in class \(s\). By
Proposition~\ref{prp-waits-across-missing-rank}, winning pairs for these
two draws must lie on opposite sides of \(p\). They both belong to class
\(s\), whose multiplicities are unchanged in \(\boldsymbol h\). This
requires two ranks of multiplicity at least \(2\) in class \(s\),
contradicting (2). Exchanging \(x\) and \(y\) proves the assertion for
\(q\), completing (1).
\end{proof}

\begin{remark}
The bound in Proposition~\ref{prp-unique-pair-ranks} (2) of at most one
rank of multiplicity \(2\) in the fixed class is attained by
\(\mathbf D_n\), whose fixed class is \(s=1\) and whose rank \(4\) has
multiplicity \(2\).
\end{remark}

\begin{proposition}[Draw-pair relation and two-pair
decomposition]\protect\hypertarget{prp-draw-pair-n-plus-one}{}\label{prp-draw-pair-n-plus-one}

With \(\boldsymbol c, x, y, s, p, q\) and the hypotheses as in
Proposition~\ref{prp-unique-pair-ranks}:

\begin{enumerate}
\def\labelenumi{\arabic{enumi}.}
\item
  The draw-pair relation of \(\boldsymbol c\) between the interchanging
  classes \(x\), \(y\) is given by \[
  \begin{aligned}
  i \to q & \text{ for every } i\equiv x\pmod3 \\
  j \to p & \text{ for every } j\equiv y\pmod3
  \end{aligned}
  \]
\item
  \(p \leftrightarrow q\) is the unique bidirectional edge in the
  draw-pair relation of \(\boldsymbol c\). \(\boldsymbol c\) has a
  unique two-pair decomposition
  \(\boldsymbol c=2\boldsymbol e_p+2\boldsymbol e_q+\boldsymbol d\), up
  to exchanging \(p\) and \(q\).
\item
  \(\boldsymbol d\) is sequentially meldable by \(n\) sequences.
  Moreover, moving a tile of \(\boldsymbol d\) from \(p\) to any rank in
  class \(x\), or from \(q\) to any rank in class \(y\), preserves its
  sequential meldability. That is,
  \(\boldsymbol d-\boldsymbol e_p+\boldsymbol e_i\) is sequentially
  meldable for every \(i\equiv x\pmod3\), and
  \(\boldsymbol d-\boldsymbol e_q+\boldsymbol e_i\) is sequentially
  meldable for every \(i\equiv y\pmod3\).
\end{enumerate}

\end{proposition}

\begin{proof}
For (1), if \(i\equiv x\pmod3\) and \(i\to j\), then
Proposition~\ref{prp-draw-pair-congruence} gives \(j\equiv y\pmod3\).
Since \(i\neq j\), we must have \(c_j\geq2\), so
Proposition~\ref{prp-unique-pair-ranks} (1) forces \(j=q\).
Full-sidedness ensures \(i\to q\). Exchanging \(x,y\) and \(p,q\) proves
the other case.

For (2), by part (1) and Proposition~\ref{prp-bidirectional-two-pair},
it suffices to rule out possible bidirectional edges in the fixed class.
Suppose \(a \to b\) and \(b \to a\) with \(a\neq b\) in the fixed class.
Then \(c_a,c_b\geq2\), contradicting
Proposition~\ref{prp-unique-pair-ranks} (2). Hence the bidirectional
edge and the two-pair decomposition are unique.

For (3), for any \(i\equiv p\pmod3\), consider a meld decomposition of
\(\boldsymbol c+\boldsymbol e_i-2\boldsymbol e_q\). Its fixed class
contains \(n\) tiles (by Proposition~\ref{prp-class-sizes-n-plus-one})
and no triplet (by Proposition~\ref{prp-unique-pair-ranks} (2)). So it
has exactly \(n\) sequences. The remaining meld is a triplet in \(x\).
By Proposition~\ref{prp-unique-pair-ranks} (1), it must be
\(\boldsymbol t_p\), as \(p\) is the only rank in class \(x\) with
multiplicity at least \(2\). Removing \(\boldsymbol t_p\) leaves
\(\boldsymbol d-\boldsymbol e_p+\boldsymbol e_i\) as a sum of \(n\)
sequences. Taking \(i=p\) proves the assertion for \(\boldsymbol d\)
itself. Symmetric arguments prove the assertions for \(q\).
\end{proof}

We now prove Theorem~\ref{thm-main}.

\begin{proof}
\textbf{The case \(n=1\).} Each congruence class contains a single rank.
The class-size bounds and congruences used in the proof of
Proposition~\ref{prp-class-sizes-n-plus-one} give
\((|\boldsymbol c|_s,|\boldsymbol c|_x,|\boldsymbol c|_y)=(1,3,3)\) or
\((3,2,2)\). Their permutations are precisely the six hands in
Theorem~\ref{thm-main}. We omit the trivial verification that each is
full-sided and live.

Henceforth let \(n\geq2\), and take \(p,q,\boldsymbol d\) from
Proposition~\ref{prp-draw-pair-n-plus-one}. Moving the tile at \(p\) to
the smallest and largest ranks in its congruence class,
Proposition~\ref{prp-draw-pair-n-plus-one} with
Proposition~\ref{prp-tile-moving} forces \(n-1\) distinct sequences to
occur in the unique sequence decomposition of \(\boldsymbol d\). The
same applies to \(q\). Since \(\boldsymbol d\) consists of only \(n\)
sequences, the two lists of required sequences must have at least
\(n-2\) sequences in common.

\textbf{The case \(s=2\).} Label \(p,q\) so that \(p\equiv1\) and
\(q\equiv0\pmod3\). The required sequence starts are \[
\begin{array}{c|ll}
 & \text{moves to the left} & \text{moves to the right} \\
\hline
p & 2,5,\ldots,p-2 & p,p+3,\ldots,3n-5 \\
q & 4,7,\ldots,q-2 & q,q+3,\ldots,3n-3
\end{array}
\] Combining the left and right entries, each row of the table contains
\(n-1\) distinct starts. The only starts that can occur in both rows are
\(4,7,\ldots,3n-5\), of which there are \(n-2\). Since at least \(n-2\)
starts must be shared, all of these occur in both rows, leaving exactly
one unshared start in each row. The left list for \(p\) and the right
list for \(q\) contain only unshared starts, so \((p-1)/3\leq1\) and
\((3n-q)/3\leq1\). Hence \(p\in\{1,4\}\) and \(q\in\{3n-3,3n\}\).

Counting the shared starts only once, the two rows together specify
\((n-1)+(n-1)-(n-2)=n\) distinct sequences. Each must occur at least
once in the sequence decomposition of \(\boldsymbol d\). But that
decomposition contains exactly \(n\) sequences in total, so every
required sequence occurs exactly once, and no other sequence occurs.
Thus the table determines \(\boldsymbol d\) completely. Summing these
sequences gives
\(\boldsymbol d=(0,1^{3n-2},0)+\boldsymbol e_p+\boldsymbol e_q\);
restoring the two pairs gives
\(\boldsymbol c=(0,1^{3n-2},0)+3\boldsymbol e_p+3\boldsymbol e_q\). The
\(4\) choices are \[
\begin{array}{c|cc}
p \operatorname{\backslash} q & 3n & 3n-3 \\
\hline
1 & \mathbf A_n & \mathbf B_n \\
4 & \mathbf B'_n & \mathbf E_n
\end{array}
\] For \(n=2\), the last entry is \(\mathbf E_2=(0,1,4,4,1,0)\).

\textbf{The case \(s=1\).} Label \(p,q\) so that \(p\equiv0\) and
\(q\equiv2\pmod3\). The required sequence starts are now \[
\begin{array}{c|ll}
 & \text{moves to the left} & \text{moves to the right} \\
\hline
p & 4,7,\ldots,p-2 & p,p+3,\ldots,3n-3 \\
q & 3,6,\ldots,q-2 & q,q+3,\ldots,3n-4
\end{array}
\] Suppose first that \(q=2\). Then \(c_2 = c_q \geq 3\), while the
other \(n-1\) ranks in class \(2\) have multiplicity \(1\) by
Proposition~\ref{prp-unique-pair-ranks} (1) and
Proposition~\ref{prp-interior-positive}. Hence \(|\boldsymbol c|_2=n+2\)
forces \(c_2=3\). If \(c_1=0\), by
Proposition~\ref{prp-unique-pair-ranks} (3), drawing rank \(1\) requires
\(\boldsymbol t_q=\boldsymbol t_2\) as a triplet in the winning
decomposition. But removing it leaves no tile of rank \(2\), leaving the
single tile of rank \(1\) impossible to form a meld or a winning pair, a
contradiction. Thus \(c_1>0\). Since \(d_1=c_1\), the sequence
\(\boldsymbol s_1\) must occur in \(\boldsymbol d\). But when \(q=2\),
the two rows of the table are disjoint and neither contains the start
\(1\). Together with \(\boldsymbol s_1\), they require \(2n-1>n\)
sequences, a contradiction.

We therefore have \(q>2\). By Proposition~\ref{prp-unique-pair-ranks}
(1) and Proposition~\ref{prp-interior-positive}, \(d_2=c_2=1\), so
\(\boldsymbol d\) contains \(\boldsymbol s_1\) or \(\boldsymbol s_2\).
Neither start appears in the table. The union of the two rows
consequently has at most \(n-1\) starts. Each row already has \(n-1\)
starts, so the rows coincide. Their only possible common starts are
divisible by \(3\); hence the left list for \(p\) and the right list for
\(q\) are empty. This forces \(p=3\) and \(q=3n-1\), and \[
\boldsymbol d=\boldsymbol s_3+\boldsymbol s_6+\cdots+\boldsymbol s_{3n-3}+\boldsymbol s_j
\qquad j\in\{1,2\}
\] Restoring \(2\boldsymbol e_3+2\boldsymbol e_{3n-1}\) gives
\(\mathbf C_n\) when \(j=1\) and \(\mathbf D_n\) when \(j=2\).

\textbf{The case \(s=0\).} A reversed argument of the previous case
gives \(\mathbf C'_n\) and \(\mathbf D'_n\).

\textbf{Verification of the eight families.} For each of the six hands
obtained in the cases \(s=2\) and \(s=1\), the displayed construction of
\(\boldsymbol d\) contains every sequence required by the two
corresponding rows. By Proposition~\ref{prp-tile-moving},
\(\boldsymbol d-\boldsymbol e_p+\boldsymbol e_i\) is sequentially
meldable for every \(i\equiv p\pmod3\). Adding \(\boldsymbol t_p\) and
the pair \(2\boldsymbol e_q\) gives a winning decomposition of
\(\boldsymbol c+\boldsymbol e_i\). Exchanging \(p,q\) verifies all draws
in the other interchanging class.

It remains to verify draws in the fixed class. Put
\(\boldsymbol h=\boldsymbol c-\boldsymbol t_p-\boldsymbol t_q\). For
\(\mathbf A_n,\mathbf B_n,\mathbf B'_n,\mathbf E_n\), the hand
\(\boldsymbol h\) consists of one tile of each rank from \(2\) to
\(3n-1\). For \(\mathbf C_n\), it consists of one tile of each rank from
\(1\) to \(3n-2\). In either case, for a fixed-class draw \(i\),
removing a pair at \(i\) leaves \(\boldsymbol h-\boldsymbol e_i\) as two
intervals of single tiles whose lengths are divisible by \(3\).
Splitting these into sequences and restoring the two triplets proves
\(i\to i\).

For \(\mathbf D_n\), the hand \(\boldsymbol h\) consists of
\(\boldsymbol s_2\) together with one tile of each rank from \(4\) to
\(3n-2\). If \(i\equiv1\pmod3\) and \(i\geq4\), then
\(\boldsymbol h-\boldsymbol e_i\) consists of \(\boldsymbol s_2\) and
two intervals of lengths divisible by \(3\), again proving \(i\to i\).
For \(i=1\), the hand \(\boldsymbol h+\boldsymbol e_1-2\boldsymbol e_4\)
consists of \(\boldsymbol s_1\) and one tile of each rank from \(5\) to
\(3n-2\), an interval of length \(3n-6\). Thus \(1\to4\). Reversing the
winning decompositions verifies \(\mathbf C'_n\) and \(\mathbf D'_n\).

This completes the proof of Theorem~\ref{thm-main}.
\end{proof}

\begin{remark}
By Proposition~\ref{prp-canonical-meld-placement}, the difference of two
full-sided \((n,n+1)\)-waits is virtually meldable if and only if they
share the same fixed class \(s\). For \(n\geq2\), this partitions the
eight families into three equivalence classes:

\begin{itemize}
\tightlist
\item
  \(s=2\): \(\mathbf A_n,\mathbf B_n,\mathbf B'_n,\mathbf E_n\);
\item
  \(s=1\): \(\mathbf C_n\) and \(\mathbf D_n\);
\item
  \(s=0\): \(\mathbf C'_n\) and \(\mathbf D'_n\).
\end{itemize}

\end{remark}

\section{\texorpdfstring{Towards the General Classification of
\(m>n+1\)}{Towards the General Classification of m\textgreater n+1}}\label{sec-amt}

With the classical configuration \(m=n+1\) completely settled, a natural
next question is the general \(m>n+1\) case. We pose the following
problem.

\begin{problem}
Classify all full-sided \((n,m)\)-waits for \(m>n+1\). Determine the
number of full-sided \((n,m)\)-waits for each \(n\) and \(m\). In
particular, solve the case of generalized Taiwanese Mahjong where
\(m=n+2\).

\end{problem}

Full-sided waits are somewhat abundant when the tile multiplicities are
large enough. In fact, we have the following result.

\begin{proposition}[Abundance of full-sided
waits]\protect\hypertarget{prp-full-sided-waits-abundant}{}\label{prp-full-sided-waits-abundant}

Let \(\boldsymbol c\) be any \(n\)-hand with
\(|\boldsymbol c| \equiv 1 \pmod 3\) and
\(\boldsymbol c \geq \boldsymbol h\), where the \(n\)-hand \[
\boldsymbol h := \begin{cases}
(2,4,2) & \text{if } n=1 \\
(2,4,6,\dots,6,4,2) & \text{if } n\geq2
\end{cases}
\] then \(p \to_{\boldsymbol c} q\) for every pair \((p,q)\) with
\(p + q \equiv -S(\boldsymbol c) \pmod 3\). In particular,
\(\boldsymbol c\) is a full-sided wait.

\end{proposition}

\begin{proof}
By hypotheses,
\(\boldsymbol d := \boldsymbol c+\boldsymbol e_p-2\boldsymbol e_q\) has
its size and rank sum congruent to \(0\) modulo \(3\), so
Proposition~\ref{prp-canonical-meld-placement} (1) gives a unique
canonical virtual meld placement of \(\boldsymbol d\). Say the placement
is \((\sigma,\tau)\). By Proposition~\ref{prp-canonical-meld-placement}
(2), to confirm that \(p \to_{\boldsymbol c} q\), it suffices to verify
\(\tau_j \geq 0\) for each \(j\). Indeed, \[
3\tau_j = d_j - \sigma_{j-2} - \sigma_{j-1} - \sigma_j \geq d_j - h_j \geq c_j - h_j - 2 \geq -2
\] since \(0 \leq \sigma_k \leq 2\) for all \(1 \leq k \leq 3n-2\) and
\(\sigma_{-1} = \sigma_0 = \sigma_{3n-1} = \sigma_{3n} = 0\). Since
\(\tau_j\) is an integer, we have \(\tau_j \geq 0\) for all \(j\), as
desired.
\end{proof}

This abundance is not a surprise, as adding melds into a wait preserves
its full-sidedness. This motivates us to reduce the problem to a more
manageable one by focusing on the \emph{minimal} full-sided waits.

\begin{definition}[Minimal full-sided
wait]\protect\hypertarget{def-minimal-full-sided-wait}{}\label{def-minimal-full-sided-wait}

A full-sided wait \(\boldsymbol c\) is \emph{minimal} if there does not
exist a full-sided wait \(\boldsymbol c'\) and a nonzero meldable hand
\(\boldsymbol d\) such that
\(\boldsymbol c=\boldsymbol c'+\boldsymbol d\).

\end{definition}

\begin{problem}
Classify minimal full-sided \((n,m)\)-waits for all \(n\) and \(m\). In
particular, determine the number of minimal full-sided \((n,m)\)-waits
for each \(n\) and \(m\).

\end{problem}

\begin{example}[Standard
Mahjong]\protect\hypertarget{exm-standard-mahjong}{}\label{exm-standard-mahjong}

Since no full-sided \((n,m)\)-wait exists for \(m\leq n\) by
Proposition~\ref{prp-no-m-le-n-full-wait}, all full-sided
\((n,n+1)\)-waits classified in Theorem~\ref{thm-main} are automatically
minimal. Hence there are \(0\) minimal full-sided \((n,m)\)-waits for
\(m\leq n\), \(8\) minimal full-sided \((n,n+1)\)-waits for \(n\geq2\),
and \(6\) minimal full-sided \((1,2)\)-waits for \(n=1\).

\end{example}

\begin{example}[Taiwanese
Mahjong]\protect\hypertarget{exm-taiwanese-mahjong}{}\label{exm-taiwanese-mahjong}

A computer search of full-sided \((3,5)\)-waits finds \(149\) solutions,
of which \(11\) are live and \(27\) are minimal, confirming
\autocite[Section 6]{cheng_mathematical_2019}'s result.

\end{example}

As a first step towards these general classification problems, we fix
\(n\) and investigate the distribution of minimal full-sided
\((n,m)\)-waits as \(m\) varies. The case \(n=1\) is relatively easy:

\begin{proposition}[]\protect\hypertarget{prp-M-m-one}{}\label{prp-M-m-one}

All minimal full-sided waits of order \(n=1\) are given by the \(6\)
hands in Theorem~\ref{thm-main} with \(m=2\).

\end{proposition}

\begin{proof}
Let \(\boldsymbol c=(c_1,c_2,c_3)\) be a full-sided wait of order \(1\).
Since for \(n=1\) any permutation of the ranks preserves melds and
full-sidedness, we may assume that rank \(2\) lies in the fixed class.
By Proposition~\ref{prp-class-size-lower-bounds} and
Proposition~\ref{prp-class-size-congruence}, \[
c_2\geq1,\qquad c_1,c_3\geq2,\qquad c_2\equiv c_1+1\equiv c_3+1\pmod3
\] Remove triplets until \(1\leq c_2\leq3\) and \(2\leq c_1,c_3\leq4\).
Removing triplets preserves the multiplicities modulo \(3\), so the
congruences above force the remaining hand to be \((3,1,3)\),
\((2,3,2)\), or \((4,2,4)\). In the last case, removing the sequence
\((1,1,1)\) leaves \((3,1,3)\). Thus, up to the initial permutation, the
remaining hand is always one of the six full-sided \((1,2)\)-waits in
Theorem~\ref{thm-main}.
\end{proof}

For the general case \(n\geq2\), we already know that \(m\) is bounded
below by \(n+1\) by Proposition~\ref{prp-no-m-le-n-full-wait}. The
following result shows that \(m\) is also bounded above.

\begin{proposition}[]\protect\hypertarget{prp-M-m-bounds}{}\label{prp-M-m-bounds}

Fix \(n\geq1\). There is a largest permitted meld count \(M(n)\) for a
minimal full-sided \((n,m)\)-wait to exist. Moreover,
\(2n\leq M(n)<8n\).

\end{proposition}

\begin{proof}
Say \(\boldsymbol c\) is a minimal full-sided \((n,m)\)-wait. In the
canonical meld placement of any winning decomposition
\(\boldsymbol c+\boldsymbol e_i-2\boldsymbol e_j\), sequences use at
most \(2+2+2=6\) tiles at any rank \(k\) and the winning pair uses at
most \(2\), so any rank with \(c_k\geq9\) would contribute a removable
triplet \(\boldsymbol t_k\) across all draws. Thus every minimal
full-sided wait satisfies \(c_k\leq8\) at each rank. Hence \(M(n)\)
exists and \(M(n)<8n\).

For a lower bound, one can check that the \(n\)-hand \((3,2^{3n-1})\) is
a minimal full-sided \((n,2n)\)-wait, giving \(M(n)\geq2n\).
\end{proof}

A precise determination of \(M(n)\) requires more work. Using exhaustive
computer search for small \(n\), we found two families of minimal
full-sided waits of even order \(n=2a+2\) and odd order \(n=2a+3\)
(\(a\geq0\)), \[
\begin{aligned}
\boldsymbol u_a &= (2,4,6,3,(5,4)^{3a},2,2) \\
\boldsymbol v_a &= (2,4,6,3,3,(5,4)^{3a+1},2,2)
\end{aligned}
\] which improve the lower bound to \[
M(n)\geq\left\lfloor\frac{9n}{2}\right\rfloor-3\quad(n\geq2)
\] We conjecture that this lower bound is in fact exact.

\begin{conjecture}[]\protect\hypertarget{cnj-upper-bound-m}{}\label{cnj-upper-bound-m}

For every \(n\geq1\), the maximum meld count \(M(n)\) of a minimal
full-sided \((n,m)\)-wait is given by \[
M(1)=2,\qquad
M(n)=\left\lfloor\frac{9n}{2}\right\rfloor-3\quad(n\geq2)
\]

\end{conjecture}

\appendix

\section*{AI Disclosure}\label{ai-disclosure}
\addcontentsline{toc}{section}{AI Disclosure}

During the preparation of this paper, the authors used OpenAI's GPT 5.6
Sol, GPT 6 Astra, and Google's Gemini 4 Argon to assist with manuscript
drafting and editing, proving auxiliary lemmas, drafting candidate Lean
4 formalizations, and exploring possible research directions. The
authors steered the process, independently verified all mathematical
arguments, and made all decisions concerning the selection,
significance, organization, and presentation of the results. The authors
fully understand the paper and take full responsibility for the content.

\ifzh
  \printbibliography[title={参考文献}]
\else
  \printbibliography
\fi

\end{document}